\documentclass{amsart}

\usepackage{amsmath,amssymb}
\usepackage{mathrsfs}
\usepackage{ifpdf}\ifpdf
\usepackage{graphicx}
\usepackage[pdftex,bookmarks=false]{hyperref}
\hypersetup{
pdfauthor={Eugenio Trucco},
pdftitle={Wandering Fatou components and algebraic Julia sets},
pdfsubject={},
pdfkeywords={Dynamical Systems}}
\else
\usepackage{hyperref}
\usepackage{graphicx}
\fi

\theoremstyle{plain}
\newtheorem{theorem}{Theorem}[section]
\newtheorem{proposition}[theorem]{Proposition}
\newtheorem{lemma}[theorem]{Lemma}
\newtheorem{corollary}[theorem]{Corollary}
\newtheorem{introtheorem}{Theorem}

\newtheorem{introcorollary}[introtheorem]{Corollary}

\theoremstyle{definition}
\newtheorem{definition}[theorem]{Definition}

\theoremstyle{remark}
\newtheorem{remark}[theorem]{Remark}

\newcommand{\A}{\mathbb A}
\newcommand{\C}{\mathbb C}

\renewcommand{\H}{\mathbb H}
\renewcommand{\L}{\mathbb L}
\newcommand{\N}{\mathbb N}

\newcommand{\Q}{\mathbb Q}
\newcommand{\R}{\mathbb R}

\newcommand{\Z}{\mathbb Z}

\newcommand{\cA}{\mathcal A}
\newcommand{\cB}{\mathcal B}
\newcommand{\cC}{\mathcal C}
\newcommand{\cD}{\mathcal D}

\newcommand{\cL}{\mathcal L}
\newcommand{\cM}{\mathcal M}
\newcommand{\cN}{\mathcal N}
\newcommand{\cO}{\mathcal O}

\newcommand{\cS}{\mathcal S}
\newcommand{\cU}{\mathcal U}
\newcommand{\cV}{\mathcal V}

\newcommand{\hD}{\mathscr D}

\newcommand{\gm}{\mathfrak m}

\newcommand{\go}{\mathfrak o}

\newcommand{\e}{\varepsilon}
\newcommand{\z}{\zeta}

\newcommand{\w}{\omega}

\newcommand{\defi}{\mathrel{\mathop:}=}

\newcommand{\tol}{\longrightarrow}

\newcommand{\barra}{\overline}
\newcommand{\htop}{h_\textnormal{top}}
\newcommand{\lec}{\preccurlyeq}
\newcommand{\sub}{\subseteq}

\newcommand{\vacio}{\varnothing}
\newcommand{\what}{\widehat}
\newcommand{\wtilde}{\widetilde}

\DeclareMathOperator{\aev}{\:a.e.}

\DeclareMathOperator{\chara}{char}

\DeclareMathOperator{\conv}{conv}
\DeclareMathOperator{\Crit}{Crit}

\DeclareMathOperator{\diam}{diam}

\DeclareMathOperator{\GO}{GO}

\DeclareMathOperator{\lcm}{lcm}

\DeclareMathOperator{\ord}{ord}

\providecommand{\abs}[1]{\left\lvert#1\right\rvert}
\providecommand{\absp}[1]{\left\lvert#1\right\rvert_p}

\providecommand{\pen}[1]{\left\lfloor#1\right\rfloor}

\newcommand{\ACP}{\A^{1,\textnormal{an}}_{\C_p}}
\newcommand{\AK}{\A^{1,\textnormal{an}}_K}

\newcommand{\CJP}{\conv(\JP)}
\DeclareMathOperator{\CritI}{Crit^{\textnormal{I}}}

\newcommand{\flau}{F(\!(\tau)\!)}
\newcommand{\flaum}{F(\!(\tau^{1/m})\!)}
\DeclareMathOperator{\lcg}{lcg}

\newcommand{\fpui}{F\langle\!\langle \tau\rangle\!\rangle}
\newcommand{\dH}{d_\H}
\newcommand{\JP}{\mathcal J_P}
\newcommand{\JPI}{\JP^{\textnormal{I}}}
\newcommand{\FP}{\mathcal F_P}	
\newcommand{\FPI}{\FP^{\textnormal{I}}}
\newcommand{\KP}{\mathcal K_P}
\newcommand{\KPI}{\KP^{\textnormal{I}}}
\newcommand{\HK}{\H_K}

\providecommand{\calg}[1]{#1^\mathrm{a}}
\providecommand{\degl}[1]{\deg_{#1}(P)}

\begin{document}

\title{Wandering Fatou Components and Algebraic Julia Sets}

\author{Eugenio Trucco}
\address{Facultad de Matem\'aticas, Pontificia Universidad Cat\'olica de Chile, Casilla 306, Correo 22, Santiago, Chile.}
\email{ejtrucco@mat.puc.cl}
\thanks{Partially supported by Research Network on Low Dimensional Dynamics, PBCT ACT-17, CONICYT, Chile and ECOS-CONICYT C07E01.}

\date{}

%\dedicatory{}

\begin{abstract}
We study the dynamics of polynomials with coefficients in a non-Archimedean field $K,$ where $K$ is a field containing a dense subset of algebraic elements
over a discrete valued field $k.$ We prove that every wandering Fatou component is contained in the basin of a periodic orbit. We obtain a complete description
of the new Julia set points that appear when passing from $K$ to the Berkovich line over $K$. We give a dynamical characterization of polynomials having algebraic
Julia sets. More precisely, we establish that a polynomial with algebraic coefficients has algebraic Julia set if every critical element is nonrecurrent.
\end{abstract}

\maketitle

\section{Introduction}

In this paper we study the dynamics of polynomials $P\colon K\to K$ where $K$ is a non-Archimedean field which is complete and algebraically closed. Moreover, we will assume
that there exists a discrete valued field $k\sub K$ such that
\[\calg k=\{z\in K\mid [k(z):k]<+\infty\}\]
form a dense subset of $K.$ Examples of such fields are the field $\C_p$ of $p$-adic numbers and the field, which we will denote by $\L,$ which is the
completion of an algebraic closure of the field of formal Laurent series with coefficients in $\C.$ Dynamics over $\C_p$ naturally arise in number theory and
dynamics over $\L$ naturally appears in the study of parameter spaces of complex rational maps~\cite{pui_kiwi}.

For complex rational maps acting on the Riemann sphere,  
Sullivan \cite{quasiconformalI_sullivan} proved, with the aid of quasi-conformal techniques, that every connected 
component of the Fatou set of a rational map $R\in\C(z)$ of degree $\ge2$ is eventually periodic (Sullivan's No Wandering Domains Theorem). 
This is no longer true
for general non-Archimedean fields. In fact, Benedetto \cite{examples_benedetto} established the existence of $p$-adic polynomials having wandering
(analytic) domains which are not attracted to a periodic orbit. This result heavily relies on the fact that over $p$-adic fields, whose residual
characteristic is $p >0$, there exists a phenomenon called~\emph{wild ramification}.

The aim of this paper is to study the interplay between algebraic and dynamical properties of points in the Julia set of a polynomial.
As a consequence, we establish that for \emph{tame} polynomials (see Definition~\ref{def:tame}), that is, for polynomials such that wild ramification does not occur, the
dynamics is free of nontrivial wandering domains (see Corollary~\ref{cor:nowandering} below).

Recent developments on the theory of iteration of rational maps over non-Archimedean fields put in evidence that
the correct space to study the action of rational maps is the \emph{Berkovich space} (e.g.
\cite{equi_bakrum,measure_chamber,theorie_favriv,tesis_rivera,tesis_aste_rivera,espace_rivera}).
The action of a polynomial $P\in K[z]$ extends to the \emph{Berkovich affine line $\AK$ associated to $K.$}
Moreover, the notions of Julia set (chaotic dynamics) and Fatou set (tame dynamics) also extend to $\AK$.
Our first main result is a complete description of the new Julia set points that appear when passing from 
$K$ to $\AK$. We will denote by $\JP$ the Julia set of $P.$

\begin{introtheorem}\label{theo:juliaberk}
Let $P\in K[z]$ be a tame polynomial of degree $d\ge2.$ Then $\JP\setminus K$ is empty or, there exist finitely
many repelling periodic orbits
$\cO_1,\dots,\cO_m \subseteq \AK \setminus K$ such that
\[\JP\setminus K=\GO(\cO_1)\sqcup\cdots\sqcup\GO(\cO_m),\]
where $\GO(\cO_j)$ denotes the grand orbit of $\cO_j.$
\end{introtheorem}

The previous theorem is first proven for polynomials in $K[z]$ with algebraic coefficients over the field $k.$ Here, we rely on our study of the interplay between
the geometry of the Julia set and the underlying algebraic structure (section~\ref{sec:polalg}). For a general tame polynomial with  coefficients in $K,$
we use a perturbation technique furnished by a key proposition (Proposition~\ref{prop:diam}) 
inspired by complex polynomial dynamics~(e.g. \cite{branner_qiuyin}).

Standard techniques (see Proposition~\ref{prop:fatoujulia}) allow us to deduce the above mentioned nonwandering result from Theorem~\ref{theo:juliaberk}.

\begin{introcorollary}\label{cor:nowandering}
Let $P\in K[z]$ be a tame polynomial of degree $\ge2.$ Then, every wandering Fatou component is in the basin of a periodic orbit.
\end{introcorollary}

Benedetto~\cite{thesis_benedetto} proved a similar result to Corollary~\ref{cor:nowandering} for rational maps with \emph{algebraic coefficients}
over $\Q_p$ with some slightly different hypothesis. 

In terms of $k,K$ and in our language, Theorem B on~\cite{wanreduc_benedetto} says that every wandering Fatou component of a rational map
with algebraic coefficients over $k$ is in the basin of a periodic orbit. Benedetto asks (question (2) at the end of the introduction 
of~\cite{wanreduc_benedetto}) if this is true for rational maps with coefficients in $K,$ assuming that the characteristic of the residual field of
$K$ is zero. Corollary~\ref{cor:nowandering} above gives an affirmative answer to the question posed by Benedetto in the case of polynomials.

After extending the notion of  \emph{algebraic degree of $x\in K$ over $k$} for arbitrary points $x \in \AK$ (see section~\ref{sec:algdeg})
we obtain the algebraic counterpart of the previous topological dynamics results.

\begin{introtheorem}\label{theo:juliaalg}
Let $P(z)\in K[z]$ be a tame polynomial of degree $\ge2$ and with algebraic coefficients over $k.$ If the algebraic degree
of every element in $\JP$ is finite then the critical elements contained in $\JP$ are not recurrent. In that case the algebraic degrees
of the elements of $\JP$ are uniformly bounded.
\end{introtheorem}

In a special type of fields, which will denote by $\L_F$, we obtain the converse of the previous theorem. Here $F$ is an algebraically closed field of characteristic 0 and $\L_F$ is the completion of an algebraic closure of the field $\flau$ of formal Laurent series with coefficients in $F$ with respect to some non-Archimedean absolute value. See subsection~\ref{subsec:lau} for definitions.

\begin{introtheorem}\label{theo:juliaalgiff}
Let $P(z)\in\L_F[z]$ be polynomial of degree $\ge2$ and with algebraic coefficients over $\flau.$ Then the algebraic degree
of every element in $\JP$ is finite if and only if the critical elements contained in $\JP$ are not recurrent. In that case, $\JP$ is contained in a
finite extension of $\flau.$
\end{introtheorem}

\subsection{Outline of the paper}

Section~\ref{sec:prelim} consists of basic definitions and facts about the Berkovich affine line and the action of polynomials on $\AK.$  

In sections~\ref{sec:dynpoints} and \ref{sec:geosec} we introduce the 
\emph{dynamical sequence} and the \emph{geometric sequence} of a polynomial. We employ these objects throughout the paper
since they organize our topological and algebraic study of the convex hull of the Julia set.

In Section~\ref{sec:algdeg} we extend  the notion of algebraic
degree to the Berkovich line and explore its basic properties as well
as the relation between the algebraic degree and the  \emph{equilibrium measure} (e.g. \cite{equi_bakrum,measure_chamber,theorie_favriv,equi_favriv})

In Section~\ref{sec:polalg}, for polynomials with algebraic
coefficients, we describe the behavior of the algebraic degree along
geometric sequences. Then we prove Corollary~\ref{cor:nowandering},
in the case of polynomials with algebraic coefficients, and finish the section with the proofs of
Theorem~\ref{theo:juliaalg} and of Theorem~\ref{theo:juliaalgiff}.

In Section \ref{sec:polele} we establish Proposition~\ref{prop:diam} which is the key to perturb polynomials (with transcendental coefficients), 
preserving the dynamics along an orbit. Then we 
prove Theorem~\ref{theo:juliaberk} in full generality and, as a consequence, we obtain Corollary~\ref{cor:nowandering}.

\section{Background}\label{sec:prelim}

\subsection{Non-Archimedean fields}

Let $K$ be a field with characteristic zero endowed with a non-Archimedean absolute value $\abs\cdot$. That is, an absolute value satisfying the
\emph{strong triangle inequality}
\[\abs{z_1+z_2}\le\max\{\abs{z_1},\abs{z_2}\}\]
for all $z_1,z_2\in K.$ Examples of such fields are the field of $p$-adic numbers $\Q_p$ and the field $\L_F$ of Puiseux series
with coefficients in $F$ that will discuss in detail in subsection~\ref{subsec:lau}. For more about non-Archimedean fields see \cite{local_cassels,padic_gouvea}. 

The set $\abs{K^\times}\defi\{\abs z\mid z\in K^\times\}$ of nonzero values of $\abs\cdot$ is a multiplicative subgroup of the positive real
numbers called the \emph{value group of $K$.} We say that the absolute value $\abs\cdot$ is \emph{discrete} if $\abs{K^\times}$ is discrete as a subset of $\R.$

We denote by $\go_K\defi\{z\in K\mid\abs z\le1\}$ the \emph{ring of integers of $K$} and by $\gm_K$ its unique maximal ideal, i.e.
$\gm_K\defi\{z\in K\mid\abs z<1\}.$ The \emph{residual field of $K$} is the quotient field $\wtilde K\defi\go_K/\gm_K.$ As we will see later, there exists a
substantial difference according the characteristic of the residual field is either 0 or $p>0.$

For $z_0\in K$ and $r>0$ we define the sets
\[B_r^+(z_0)\defi\{z\in K\mid\abs{z-z_0}\le r\}\quad\textnormal{   and   }\quad B_r(z_0)\defi\{z\in K\mid\abs{z-z_0}<r\}.\]
If $r$ belongs to the value group of $K,$ then the sets defined above are different and we say that $B_r^+(z_0)$ (resp. $B_r(z_0)$) is a \emph{closed ball}
(resp. \emph{open ball}). If $r$ is not in the value group of $K,$ then the sets $B_r(z_0)$ and $B_r^+(z_0)$ coincide and we say that
$B_r(z_0)=B_r^+(z_0)$ is an \emph{irrational ball.} Despite these names, every ball is open and closed in the metric topology induced in $K$ by the absolute
value $\abs\cdot.$

\subsection{Balls and polynomial}

Given a nonconstant polynomial $P$ with coefficients in $K,$ define the \emph{local degree} of $P$ at $z_0\in K$ as the largest integer
$\degl{z_0}\ge1$ such that $(z-z_0)^{\degl{z_0}}$ divides $P(z)-P(z_0)$ in the ring $K[z].$ If $\degl{z_0}>1$ we say that
$z_0$ is a \emph{critical point of $P$} with \emph{multiplicity} $\degl{z_0}-1.$ We denote by $\CritI(P)$ the subset of $K$ formed by the critical
points of $P.$

The image of a ball $B\sub K,$ under the action of $P,$ is a ball, of the same type than $B$, and there exists an integer larger than 1,
denoted by $\degl B,$ and called \emph{degree of $P$ at $B,$} such that
\[\degl B=\sum_{\{z\in B\mid P(z)=z'\}}\degl z\]
for all $z'\in B$ (e.g. see section~2 in \cite{tesis_aste_rivera}).

Moreover, the preimage of a ball $B$ is a finite union of pairwise disjoint balls $B_1,\dots,B_m$ of the same type than $B$ and
\[\sum_{j=1}^m\degl{B_j}=\deg(P).\]

\subsection{Berkovich affine line}\label{subsec:berk}

We will need only basic facts about the structure of the Berkovich affine line and its topology. For more details see
\cite{valuative_favjon,tesis_aste_rivera,espace_rivera,notes_rivera}, for the original construction of
V. G. Berkovich, see \cite{spectral_berkovich}.

We identify the Berkovich line with an appropriate quotient of the set $\cS_K$ of all the strictly decreasing sequences of closed balls of $K.$ This construction is a slight modification of the given one in \cite{notes_rivera}.

On the set $\cS_K$ we define the equivalence relation $\sim$ given by: $(B_j)\sim(B_j')$ if for all $n\in\N,$ the sequence $(B_j)$ (resp. $(B_j')$)
is eventually contained in $B_n'$ (resp. $B_n$).

The \emph{Berkovich analytic space associated to the affine line over $K$} (for short, the Berkovich line) denoted by $\AK,$ is (as a set) the
quotient 
$\cS_K/\sim.$

If the sequence $(B_j)$ is equivalent to $(B_j'),$ then $\cap B_j=\cap B_j'.$ Note that the field $K$ is not spherically complete, that is,
there exist decreasing sequences of closed balls having empty intersection. However, consider $(B_j)\in\cS_K$ such that $B=\cap B_j$ is not empty.
Then $B$ is a closed ball, an irrational ball or a point of $K.$ Moreover, the intersection $B$ determines completely the equivalence class of
$(B_j)$.
In this case, we denote the equivalence class of $(B_j)$ by $x_B$ and we will say that $x_B$ in $\AK,$ is \emph{the point associated to $B$} and that $B$ is \emph{the 
ball associated to $x_B.$}

The elements of the Berkovich line $\AK$ are classified in the following four types:
\begin{enumerate}
\item \emph{Type I} or \emph{classical points,} corresponding to the equivalence classes of sequences whose intersection is a point in $K.$ We
identify $K$ with these elements of $\AK.$

\item \emph{Type II} or \emph{rational points,} corresponding to elements $x_B$ where $B$ is a closed ball.

\item \emph{Type III} or \emph{irrational points} that is, the points of the form $x_B$ where $B$ an irrational ball of $K.$

\item \emph{Type IV} or \emph{singular points,} corresponding to the equivalence classes of decreasing sequences of closed balls with empty
intersection.
\end{enumerate}

The inclusion between the balls of $K$ induces a partial order, denoted by $\lec,$ in $\AK.$ If $x\in\AK$ (resp $x'\in\AK$) is the equivalence
class of $(B_j)$ (resp. $(B_j')$) we say that $x\lec x'$ if, for each $n\in\N,$ the sequence $(B_j)$ is eventually contained in $B_n'.$ We say that $x\prec y$ if $x\lec y$ and $x\ne y.$ In the case that $x_B$ and $x_{B'}$ are nonsingular elements, we have that $x_B\lec x_{B'}$ if and only if $B$ is contained in $B'.$

For all $x\in\AK$ denote the set of elements larger than $x$ by
\[[x,\infty[\defi\{w\in\AK\mid x\lec w\}.\]
Observe that $[x,\infty[$ is isomorphic, as an ordered set, to $[0,+\infty[\sub\R.$

Given two points $x,y$ in the Berkovich line $\AK$ we have that
\[[x,\infty[\,\cap\,[y,\infty[\,=[x\vee y,\infty[\]
where $x\vee y$ is the smallest element larger than $x$ and $y.$ If $x$ is different than $y$ then the element $x\vee y$ is a type II point.

Given two elements $x,y\in\AK$ let
\[[x,y]\defi\{w\in\AK\mid x\lec w\lec x\vee y\}\cup\{w\in\AK\mid y\lec w\lec x\vee y\}.\]
The sets $]x,y],[x,y[$ and $]x,y[$ are defined in the obvious way.

For $x$ in the Berkovich line, the \emph{diameter of $x$} is
\[\diam(x)\defi\lim_{j\to\infty}\diam(B_j),\]
where $(B_j)$ is a representative of $x.$ For $x_B,$ a nonsingular element, the diameter of $x_B$ coincides with the diameter (radius) of the ball $B.$

In order to endow the Berkovich affine line with a topology, we define an \emph{open ball of $\AK$} and a \emph{closed ball of $\AK$} by
\begin{align*}
\cB(a,r)&=\{x\in\AK\mid\diam(a\vee x)<r\}\\
\cB^+(a,r)&=\{x\in\AK\mid\diam(a\vee x)\le r\}
\end{align*}
respectively, where $a\in K$ and $r>0.$

The \emph{weak topology} on the Berkovich line is the smallest topology containing all the open balls and the complements of closed balls of $\AK$.

If $B=B_r^+(a)\sub K$ we have that the closure $\barra B$ of $B$ in $\AK$ is $\cB^+(a,r).$ The boundary of $\cB^+(a,r)$ is $\{x_B\},$ although we will often
abuse of notation and write simple $\partial\cB^+(a,r)=x_B.$

For all $x\ne y\in\AK$ there exists an order preserving bijection between $[x,\infty[$ and an interval of $\R.$ Moreover, there exists an isomorphism
between $[x,x\vee y]$ and a closed interval of $\R.$ Hence, following Definition~3.5 in \cite{valuative_favjon} $\AK$ is an \emph{unrooted nonmetric
tree.}

Let $B=B_r^+(a)$ be a closed ball of $K$ and consider $x_B\in\AK$ the type II point associated to $B.$ We say that two elements $x,y\prec x_B$ are in the
same \emph{direction at $x_B$} if $x\vee y\prec x_B.$ Given $x\prec x_B,$ the set of elements in the same direction as $x$ at $x_B$ is the open
ball $\cB(z,r)$ of the Berkovich line, where $z\in B$ is such that $z\prec x.$

The \emph{tangent space at $x_B$,} denoted by $T_{x_B},$ is the set of all the directions at $x_B.$

After an affine change of coordinates $h,$ such that $h(B_1^+(0))=B$ we can identify the directions in $T_{x_B}$ with the directions at the point
associated to the ball $B_1^+(0)$ and these directions can be naturally identified with the residual field of $K.$

We say that a set $X$ of $\AK$ is \emph{convex,} if for all $x,y\in X$ we have that $[x,y]$ is contained in $X.$ For $X$ a subset of $\AK$ we
define the \emph{convex hull of $X$} to be the set
\[\conv(X)=\bigcup_{x,y\in X}[x,y].\]
A convex subset of $\AK$ is always a connected set.

\subsection{The action of a polynomial over the Berkovich line}\label{subsec:actberk}

The action of a nonconstant polynomial $P$ with coefficients in $K$ has a unique continuous extension to $\AK,$ which we also denote by $P.$
More precisely, if $x$ is the equivalence class of $(B_j),$ then $P(x)$ is defined as the equivalence class of the sequence $(P(B_j)).$ If $x=x_B$ is
a non
singular element, we have that $P(x_B)=x_{P(B)}.$

The map $P\colon\AK \to \AK$ is increasing, open and preserves the type of the points. For all $x\in\AK$ the set of preimages of $x$ under $P$ is
finite. The image of a ball $\cB$ of $\AK$ is a ball of the same type, and its preimage is a finite union of pairwise disjoint balls
of the same type than $\cB.$

To extend the notion of local degree of $P$ to $\AK$ let
\[\degl x\defi\lim_{j\to\infty}\degl{B_j},\]
where $(B_j)$ is a representative of the class $x.$ We have that $\degl{x_B}=\degl B,$ for all non
singular elements $x_B\in\AK.$ 

\begin{remark}\label{rem:sumad}
Given $x\in\AK$ with preimages $x_1,\dots,x_m,$ we have that
\[\sum_{j=1}^m\degl{x_j}=\deg(P).\]
\end{remark}

We say that $x\in\AK$ is a \emph{critical element of $P$} if $\degl x\ge2.$ The structure of the critical set 
\[\Crit(P)\defi\{x\in\AK\mid\degl x\ge2\}\]
depends strongly on the characteristic of the residual field $\wtilde K,$ as we will see in subsection~\ref{subsec:critset}. 

Let $x\in\AK$ be a type II point. Given a direction $\cD$ in $T_x,$ that is, an open ball $\cD$ of $\AK$ such that $\partial\cD=x,$ we have that $P(\cD)$ is a direction in $T_{P(x)}.$ Hence, the action
of $P$ in the Berkovich line induces a map $T_xP\colon T_x\to T_{P(x)}$ between the tangent spaces at $x$ and $P(x).$ After affine changes of
coordinates $h_1,h_2$ such that $h_2(P(B))=\go$ and 
$h_1(\go)=B$ the map $T_xP$ coincides with the reduction of $P$ to the residual field $\wtilde K.$ Hence, $T_xP$ is a polynomial map in $\wtilde K[z]$ of 
degree lower than or equal to $\deg(P).$ 

For further reference we establish a relation between the local degree of $P$ at a type II point $x\in\AK$ and the degree of $T_xP.$

\begin{remark}\label{rem:gradotangentetilde}
Let $P\in K[z]$ be polynomial of degree $\ge2$ and let $x\in\AK$ be a type II point. If $\z'\in\wtilde K$ then
\[\degl x=\deg(T_xP)=\sum_{\{\z\in \wtilde K\mid T_xP(\z)=\z'\}}\deg_{\z}(T_xP).\]
\end{remark}

\subsection{The Hyperbolic Space}

We denote by $\HK$ the \emph{hyperbolic space of $K$}, that is, the set of nonclassical elements in the Berkovich Line. This set has a tree
structure induced by the structure of $\AK.$

Over $\HK$ we can define the \emph{hyperbolic distance,}
\[\dH(x,y)=2\log\diam(x\vee y)-\log\diam(x)-\log\diam(y),\]
which is compatible with the tree structure of $\HK.$ More precisely, the set $\HK$ with the hyperbolic distance is a $\R$-tree. That is, for
all $x,y\in\HK$ the length of the segment $[x,y],$ which is a geodesic segment, coincides with the hyperbolic distance between $x$ and $y.$

For further reference we state, without proof, the following straightforward fact

\begin{lemma}\label{lemma:distlin}
Let $w,y,x$ in $\AK.$ We have that $\dH(x,y)=\dH(x,w)+\dH(w,y)$ if and only if $w$ belongs to $[x,y].$
\end{lemma}

The hyperbolic distance behaves nicely under the action of a polynomial. More precisely we have the following lemma which is a restatement of
Corollary~4.8~of~\cite{espace_rivera}.

\begin{lemma}\label{lemma:distgrado}
Let $P\in K[z]$ and consider $x\prec x'\in\AK.$ Suppose that $\degl y=\lambda$ for all $y\in\,]x,x'[.$ Then
\[\dH(P(x),P(x'))=\lambda\cdot\dH(x,x').\]
\end{lemma}

The metric topology induced in $\HK$ by the hyperbolic distance is called \emph{strong topology.} Every open set $X\sub\HK$ for the topology
induced in $\HK$ by the weak topology is an open set for the strong topology. Moreover, $(\HK,\dH)$ is a complete metric space.

Our default topology will always be the weak topology in $\AK$ and $\HK.$

\subsection{The Critical Set of $P$}\label{subsec:critset}

Let $P\in K[z]$ be a polynomial of degree $\ge2.$ The structure of the critical set of $P$ depends strongly on the characteristic of the
residual field of $K,$ as we will see in the following propositions. We will first assume that $\chara(\wtilde K)=0.$

\begin{proposition}
Let $P\in K[z]$ be a polynomial of degree $\ge 2$ and let $B\sub K$ be a ball. If the characteristic of $\wtilde K$ is zero, then
\[\degl B=1+\deg_B(P')\]
\end{proposition}
\begin{proof}
After an affine change of coordinates we can suppose that $B$ and $P(B)$ contain 0. Since $\chara(\wtilde K)=0$ we have that $\abs{n}=1$ for all $n\in\N.$ Therefore, the Newton polygon of $P'$ is a translation of the Newton polygon of $P-P(0).$ Hence, the number of zeros of $P'$ in $B$ is the number of zeros of $P$ in $B$ minus 1.
\end{proof}

\begin{remark}\label{rem:RH}
Note that if $\chara(\wtilde K)=0$ and $B\sub K$ is a ball, then the following holds
\begin{align*}
\degl B-1&=\sum_{z\in B}\Big(\degl z-1\Big)\\
&=\sum_{w\in B\cap\CritI(P)}\Big(\degl w-1\Big).
\end{align*}
That is, the degree of $P$ at the ball $B$ is determined by the critical points of $P$ contained in $B.$ We will refer to the identity above as the \emph{Riemann-Hurwitz formula}.
\end{remark}

From above, we have that in the case of $\chara(\wtilde K)=0,$ the set $\Crit(P)$ coincides with the set
\[\bigcup_{w\in\CritI(P)}[w,\infty).\]
Therefore $\Crit(P)$ is a finite subtree of $\AK.$ That is, $\Crit(P)$ has finitely many vertices and finitely many edges. There is one distinguished edge of the form
$[x,\infty[.$ The other edges are closed segments. In particular, we have that the local degree at a singular element is always 1.

The situation in the case of $\chara(\wtilde K)=p>0$ is, in general, completely different.

\begin{proposition}
Let $P\in K[z]$ be a polynomial with degree $\ge 2$ and let $B\sub K$ be a ball such that $m=\degl B$. If the characteristic of $\wtilde K$ is $p>0$ and $(p,m)=1$ then
\[\degl B=1+\deg_B(P')\]
\end{proposition}
\begin{proof}
After an affine change of coordinates we can suppose that $0$ belongs to $B$ and $P(B).$ Since the local degree is $m,$ we have that $(m,\log(\abs{a_m}))$ is a
vertex of the Newton polygon of $P.$ Since $(p,m)=1$ it follows that
\[(m-1,\log(\abs{ma_m}))=(m-1,\log(\abs{a_m}))\]
is a vertex of the newton polygon of $P'.$ In the Newton polygon of $P'$ the slope before $m-1$ can only increase and the slope after $m-1$ can only decrease
with respect to to the slopes before and after $m$ in the Newton polygon of $P.$ Hence, the number of zeros of $P'$ in $B$ is the number of zeros of $P$ in $B$ minus 1.
\end{proof}

\begin{proposition}\label{prop:pnofinito}
Let $P\in K[z]$ be a polynomial with degree $\ge 2$ and let $B\sub K$ be a closed or irrational ball such that $m=\degl B$. If the characteristic of
$\wtilde K$ is $p>0$ and $m=p^rn$ with $(p,n)=1$ and $r>1,$ then $\Crit(P)\cap\HK$ has nonempty interior with respect to the strong topology.
In particular, $\Crit(P)$ is not a finite tree.
\end{proposition}
\begin{proof}
Let $x_B$ the point associated to the ball $B.$ Since $\degl{x_B}=p^rm$ there exists a type II point $x_B\prec x$ such that $\degl x$ is also $p^rn$
and $\degl{\cD}=p^rn,$ where ${\cD}$ is the direction at $x$ that contains $x_B.$ If we consider the action between $T_x$ and $T_{P(x)},$ we have that
$p^rn=\deg(T_xP)=\deg_{\cD}(T_xP).$ After affine changes of coordinates we can suppose that $x$ and $P(x)$ are the point associated to the ball $B_1^+(0)$ and $0\in \cD=P(\cD).$
It follows that $T_xP,$ which is a polynomial of degree $p^rn,$ has a fixed point with local degree $p^rn.$ Hence  
\[T_xP(\z)=\zeta^{p^rn}=(\z^p)^{p^{r-1}n}.\]
By Lemma 10.1 in \cite{points_rivera} we have that $P$ coincide with $z^{p^rn}$ in a strong neighborhood $\cU$ of $x.$ Since $x$ is a \emph{inseparable fixed point} (Definition~5.4~in~\cite{points_rivera}) for $z^{p^rn}$
we can use Proposition 10.2 in \cite{points_rivera} to conclude the existence of a strong neighborhood $\cV$ of $x$ such that $\deg_y(z^{p^rn})>p$ for all $y\in\cV.$ Then
\[\degl y=\deg_y(z^{p^rn})>p\]
for all $y\in\cU\cap\cV.$ Therefore $\Crit(P)$ has nonempty interior with respect to the strong topology. In this case $\Crit(P)$ is not a finite tree.    
\end{proof}

The following definition is motivated by the previous propositions

\begin{definition}\label{def:tame}
We say that a nonsimple polynomial $P\in K[z]$ is \emph{tame} if the critical set of $P$ is a finite tree. 
\end{definition}

For instance, if $\chara(\wtilde K)=0$ then any polynomial is tame. If the residual characteristic of $K$ is $p>0,$ then any polynomial with degree $d<p$ is tame. 

If the residual characteristic of $K$ is $p>0$ we have, by Proposition~\ref{prop:pnofinito}, that $P$ is a tame polynomial if and only if $P$ is a nonsimple polynomial and the set
\[\{\degl x\mid x\in\AK\}\]
does not contains multiples of $p.$

If $P\in K[z]$ is a tame polynomial we have that $\Crit(P)$ coincides with the set
\[\bigcup_{w\in\CritI(P)}[w,\infty).\]
Moreover, the Riemann-Hurwitz formula (see Remark~\ref{rem:RH}) is valid. Our main results will be on tame polynomials.

\subsection{Julia and Fatou sets in the Berkovich line}

In analogy to complex polynomial dynamics the \emph{filled Julia set of $P$} is defined by
\[\KP\defi\{x\in\AK\mid(P^n(x))\,\textnormal{ is precompact}\}.\]
The filled Julia set of $P$ is always nonempty, since it contains the classical periodic points of $P.$

We define the \emph{Julia set of P,} denoted by $\JP,$ as the boundary of the filled Julia set of $P,$ that is, 
$\JP=\partial{\KP}.$ An equivalent definition, which will be useful, is the following: a point $x\in\AK$ belongs to $\JP$ if for every open
neighborhood $V$ of $x,$ we have that
\[\AK\setminus\bigcup_{n\ge0}P^n(V),\]
has at most one element.

The Julia set is a compact, totally invariant (i.e. $P(\JP)=\JP=P^{-1}(\JP)$) and nonempty set. Moreover, for all $n\in\N$ we have
$\JP=\mathcal J_{P^n}.$ Furthermore, it can be characterized as the smallest compact set totally invariant by the action of $P.$

\emph{The Fatou set of $P,$} denoted by $\FP,$ is defined as the complement of the Julia set of $P.$ This is a nonempty open set. 
We say that $\AK\setminus\KP$ is the \emph{basin of attraction of $\infty.$} Note that the basin of attraction of $\infty$ is a convex set, and 
therefore a connected set. Moreover, it is a Fatou component.

The \emph{classical filled Julia set of $P$,} denoted by $\KPI,$ is defined as $\KP\cap K.$ We define the classical Julia set of $P$ as
$\JPI\defi\JP\cap K.$ The \emph{classical Fatou set} $\FPI$ is the intersection between $\FP$ and $K.$ These definitions, of classical Fatou and
Julia sets, agree with the ones given by Hsia~\cite{weak_hsia,closure_hsia}.

Consider $x$ in $\AK$ a periodic element of period $q.$ In the case that $x$ belongs to $K,$ we say
that $x$ is \emph{attracting, neutral} or \emph{repelling} according $\abs{(P^q)'(x)}<1,$ $\abs{(P^q)'(x)}=1$ or $\abs{(P^q)'(x)}>1,$ respectively.
If $x$ belong to $\AK\setminus K$ we say that $x$ is \emph{neutral} or \emph{repelling} if $\deg_x(P^q)=1$ or $\deg_x(P^q)\ge2.$

A periodic point $x\in\AK$ of $P$ belongs to the Julia set of $P$ if and only if it is a repelling periodic point.

We will use the following proposition which is proved in section 5 of \cite{espace_rivera}.

\begin{proposition}\label{prop:critratintro}
Let $P\in K[x]$ be a polynomial of degree $\ge2$ and let $x$ be in the Julia set of $P.$ If $x$ is a periodic critical element then $x$ is a type II
point.
\end{proposition}

We say that a polynomial $P\in K[z]$ of degree $d\ge2$ is \emph{simple} if there exists a fixed point $x\in\HK$ with $\degl x=d.$

The simplest Julia set consists of a unique type II point in $\HK$ which is fixed under $P.$ In fact, the polynomials with a unique type II point as
Julia set are precisely the simple polynomials. Moreover, a tame polynomial $P$ is simple, if and only if all the classical critical points of $P$ belong
to $\KPI$ (see Corollary 2.11 in \cite{pui_kiwi}, the proof of that corollary is valid for tame polynomials).

From subsection~\ref{subsec:critset} we have that the critical set of $P$ always contains infinitely many elements. Nevertheless, if $P$ is a tame polynomial,
there are only finitely many critical elements of $P$ contained its Julia set, that is one of the important properties of tame polynomials.

\begin{proposition}\label{prop:finitecritical}
Let $P\in K[z]$ be a tame polynomial of degree $d\ge2.$ Then $\JP$ contains at most $d-2$ critical points of $P.$
\end{proposition}

In order to give the proof of the previous proposition we need the following lemma.

\begin{lemma}\label{lemma:segmentfatou}
Let $P\in K[z]$ be a nonsimple polynomial of degree $\ge2$ and consider $x\in\JP.$ Then $A(x)=\{y\in\AK\mid y\prec x\}\sub\FP\cap\KP$ and
$]x,\infty[$ is contained in the basin of $\infty.$
\end{lemma}
\begin{proof}
If $x$ is a classical or a singular point then $\vacio=A(x)\sub\FP.$ Suppose that $x=x_B\in\JP\cap\HK$ is a nonsingular point and
consider $y\in A(x).$ Then, for an open ball $\cB$ such that $y\in\cB\sub\barra B$ and given $\cD$ a ball of the Berkovich line such that 
$\JP\sub\cD$ we have that
\[P^n(\cB)\sub P^n(\barra B)\sub\cD\]
for all $n\ge1.$ That is $y$ belongs to $\FP\cap\KP.$ Therefore $A(x)\sub\FP.$

If there exists $y\in\JP$ such that $x\prec y,$ then $x\in A(y)\sub\FP,$ which is impossible. It follows that $]x,\infty[\,\sub\FP.$
\end{proof}

\begin{proof}[Proof of Proposition~\ref{prop:finitecritical}]
Recall that the critical points of $P$ belong to
\[\bigcup_{w\in\CritI(P)}[w,\infty[.\]
For each critical element $c\in\JP,$ we can choose $w_c\in\CritI(P)$ such that $c$ belongs to $[w_c,\infty[.$ In view of
Lemma~\ref{lemma:segmentfatou} we have that $[w_c,c[$ and $]c,\infty[$ are contained in the Fatou set of $P.$ Then, the map $c\mapsto w_c$ is
injective. It follows that $\JP$ contains at most $d-1$ critical elements, since $\CritI(P)$ contains at most $d-1$ elements.

Seeking a contradiction suppose that $J_P\cap\Crit(P)$ contains $d-1$ elements, it follows that $\CritI(P)$ is a subset of $\KPI.$ Following
Corollary~2.11 in \cite{pui_kiwi} we have that $P$ is a simple polynomial and the Proposition follows.
\end{proof}

The following proposition shows that the existence of wandering Fatou components is equivalent with the existence of nonpreperiodic
points in $\JP\cap\HK.$

\begin{proposition}\label{prop:fatoujulia}
Let $P\in K[z]$ be a nonsimple polynomial of degree $\ge2.$ There exists a wandering component of $\FP$ which is not in the basin of a
periodic orbit if and only if there exists a nonpreperiodic point of type II or III in $\JP.$
\end{proposition}
\begin{proof}
First note that the Fatou components different from the basin of $\infty$ are open balls of the Berkovich affine line. 

We proceed by contradiction. Let $\cB$ be a wandering Fatou component which is not in the basin of a periodic point and suppose that 
$x=\partial\cB$ is a preperiodic point. It follows that $\cB$ belongs to the basin of the orbit of $x,$ which is a contradiction.  

Conversely, if $x\in\JP\cap\HK$ is a nonsingular point which is nonpreperiodic, we have that for each $a\in K$ with $a\lec x$ the open ball
$\cB(a,\diam(x))\sub\FP\cap\KP$ is a wandering Fatou component which is not in the basin of a periodic orbit.
\end{proof}

For more results about Julia and Fatou set for rational maps see \cite{theorie_favriv,espace_rivera}.

\subsection{Measure on the Berkovich affine line}

Given a polynomial $P\in K[z]$ of degree $\ge2,$ Favre and Rivera \cite{theorie_favriv} construct an ergodic probability measure, defined
on the Borel sets of $\AK.$ See also \cite{equi_bakrum,measure_chamber,equi_favriv}. This measure is denoted by $\rho_P$ and called \emph{the
equilibrium measure of $P$}. The measure $\rho_P$ is characterized by the
following property: if $\cB$ is a ball of $\AK$ then
\[\rho_P(\cB)=\frac{\degl{\cB}}{\deg(P)}\rho_P(P(\cB)).\]

The equilibrium measure of $P$ is supported on the Julia set of $P$ and is an atom free measure for all $P$ which are not simple. Moreover, for any
open
set $V$ such that $\JP\cap V\ne\vacio$ we have that $\rho_P(V)>0.$

\section{Dynamical Points.}\label{sec:dynpoints}

Consider a nonsimple polynomial $P\in K[z]$ of degree $\ge2.$ To establish properties about $\JP$ we study the action of $P$ in the convex
hull of its Julia set, that is
\[\CJP=\bigcup_{x,y\in\JP}[x,y].\]

For each $x\in\JP$ we will construct a decreasing sequence $(L_n(x))\sub\CJP$ of type II points having $x$ as its limit. This sequence is 
dynamically defined, therefore every dynamical property of $x$ can be obtained from the properties of the sequence $(L_n(x)).$ Compare with the
\emph{lemniscates} in \cite{compacite_bezivin} and the \emph{dynamical ends} in \cite{pui_kiwi}.

At the end of the section we will introduce the concept of \emph{good starting level}, which will be useful to compare the distances between the
points in $(L_n(x))$ and $(L_n(P(x))).$

\vspace{.4cm}

From Proposition 6.7 in \cite{tesis_rivera} we know that
\[r_P=\max\{\abs{z_0-z_1}\mid z_0,z_1\in\KPI\}\]
belongs to the value group of $K.$ Thus, the closed ball $D_0=B_{r_P}^+(z)\sub K,$ where $z$ is any periodic point of $P,$ is the smallest ball of $K$ containing $\KPI,$ and therefore
$\barra{D_0}$ is the smallest ball of $\AK$ containing $\KP.$ In particular we have that the Julia set of $P$ is contained in $\barra{D_0}.$

\begin{lemma}
Let $P\in K[z]$ be a nonsimple polynomial of degree $\ge2.$ Then $\barra{D_0}$ is the smallest ball of $\AK$ containing $\JP.$ 
\end{lemma}
\begin{proof}
We proceed by contradiction. Suppose that there exists a ball $\cB$ such that $\JP\sub\cB\subsetneqq\barra{D_0}.$ From
Lemma~\ref{lemma:segmentfatou} we have that
\[\KP=\bigsqcup_{x\in\JP}\{y\in\AK\mid y\lec x\}.\]
Hence, $\KPI=\KP\cap K\sub\cB\cap K,$ which contradicts that $\diam(\KPI)=\diam(D_0).$
\end{proof}

\begin{definition}
The \emph{level 0 dynamical point of $P$}, denoted by $L_0,$ is defined as the point associated to the ball $D_0,$ that is, $L_0\defi
x_{D_0}=\partial\barra{D_0}.$
\end{definition}

\begin{definition}
For each $n\in\N$ the \emph{level $n$ dynamical set of $P$} is defined as $\cL_n\defi P^{-n}(L_0).$ We say that a element
$L_n$ of $\cL_n$ is a \emph{level $n$ dynamical point of $P.$}
\end{definition}

From the definition we have that $L_0$ is a type II point and that $x\lec L_0$ for all $x\in\JP.$ Moreover, $L_0$ is the smallest
element in $\AK$ with this property.

\begin{proposition}\label{prop:l0}
Let $P\in K[z]$ be a nonsimple polynomial of degree $d\ge2.$ Then the following statements hold:
\end{proposition}
\begin{enumerate}
\item $\{L_0\}=P^{-1}(P(L_0))$\label{item:l0ppl0propl0}

\item $L_0\prec P(L_0).$\label{item:l0precpl0propl0}

\item $\diam(P^n(L_0))\tol+\infty$ as $n\to+\infty.$\label{item:diaminfipropl0}

\item $P^{-1}(L_0)$ has at least two elements. Moreover, the elements of $P^{-1}(L_0)$ are pairwise incomparable with respect to $\lec.$\label{item:pmenospropl0}

\item $P^{-1}(L_0)$ contains points in at least two directions in $T_{L_0}.$\label{item:l0branchpropl0}
\end{enumerate}
\begin{proof}
To prove the first statement, note that, the Julia set of $P$ is forward invariant, therefore
\[\JP=P(\JP)\sub P(\barra D_0)=\barra{P(D_0)}.\]
By definition of $L_0$, we have $L_0\lec P(L_0).$ Now seeking a contradiction, suppose that there exists $x_B$ in $P^{-1}(P(L_0))$ 
different than $L_0.$ Since $\JP\sub P(\barra D_0)$ we have that $\barra B\cap\JP\ne\vacio.$ Hence, $x_B$ is comparable to $L_0.$ If $x_B\prec L_0$ it follows that
\[P(L_0)=P(x_B)\prec P(L_0),\]
which is impossible. Analogously is we suppose that $L_0\prec x_B.$ Therefore $x_B=L_0,$ which is a contradiction. Hence, we have proved that $\{L_0\}=P^{-1}(P(L_0)).$

To prove (\ref{item:l0precpl0propl0}) note that $L_0\lec P(L_0)$ and $\{L_0\}=P^{-1}(P(L_0)).$ Using Remark~\ref{rem:sumad} we obtain
that $\degl{L_0}=d,$ because $L_0$ is the unique preimage of $P(L_0).$ Hence $L_0\prec P(L_0),$ since $P$ is not a simple polynomial.

In order to show (iii) let $0<a=\dH(L_0,P(L_0)).$ We proved that $\degl{L_0}=d,$ therefore $\degl{P^n(L_0)}=d$ for
all $n\in\N.$ Following Lemma~\ref{lemma:distgrado}, we obtain that $\dH(P^{n-1}(L_0),P^{n}(L_0))=a\cdot d^{n-1}.$ Hence, $\dH(L_0,P^{n}(L_0))>a\cdot
d^{n-1},$
because all the iterates of $L_0$ belong to the segment $[L_0,\infty[.$ Using the definition of the hyperbolic distance we obtain
\[\diam(P^{n}(L_0))>\diam(L_0)\exp(a\cdot d^{n-1}).\]
Then, $\diam(P^{n}(L_0))\to+\infty.$

To prove the first statement in (\ref{item:pmenospropl0}) we proceed by contradiction. Suppose that $P^{-1}(L_0)$ has exactly one element $x_B.$ It follows that the Julia
set
of $P$ is contained in $\barra B,$ since $\JP$ is totally invariant. By the definition of $L_0$ we have that $L_0\lec x_B.$ By monotonicity of $P,$ we
obtain that $P(L_0)\lec P(x_B)=L_0,$ which is a contradiction with (\ref{item:l0precpl0propl0}). That is, $P^{-1}(L_0)$ contains at least two element.

Suppose now that there exist $x_1,x_2\in P^{-1}(L_2)$ with $x_1\prec x_2.$ It follows that
\[L_0=P(x_1)\prec P(x_2)=L_0,\]
which is impossible. Therefore the elements in $P^{-1}(L_0)$ are pairwise  incomparable.

To show (v) note that if $P^{-1}(L_0)=\{y_1,\dots,y_m\}$ we have that $x\lec y_1\vee\cdots\vee y_m$ for all $x\in\JP,$ since,
\[\JP=\bigsqcup_{1\le j\le m}\{x\in\JP\mid x\lec y_j\}.\]
In particular, if $P^{-1}(L_0)$ is contained in a direction $\cD\in T_{L_0}$ we have that
\[y_1\vee\cdots\vee y_m\prec L_0,\]
which contradicts the definition of $L_0.$ Now (v) follows.
\end{proof}

From (\ref{item:l0ppl0propl0}), (\ref{item:l0precpl0propl0}), (\ref{item:pmenospropl0}) of Proposition~\ref{prop:l0} we have that each level
$n\ge1$ dynamical point is strictly smaller than exactly one level $n-1$ dynamical point.

\begin{definition}
A \emph{dynamical sequence} is a decreasing sequence $(L_n)_{n\ge0}$ of dynamical points such that $L_0$ is the level 0 point and $L_n\in\cL_n$ for all
$n\ge1.$
\end{definition}

\begin{proposition}\label{prop:juliasequence}
Let $P\in K[z]$ be a nonsimple polynomial of degree $d\ge2.$ Then
\[\JP=\{\lim L_n\mid(L_n)\textnormal{ is a dynamical sequence of }P\}.\]
\end{proposition}
\begin{proof}
Let $(L_n)$ be a dynamical sequence of $P$ and $x=\lim L_n.$ For all $n\ge0$ we have that $x\prec L_n,$ therefore $P^n(x)\prec L_0.$ Hence $x$
belongs to $\KP.$ The dynamical points do not belong to $\KP$ (see Proposition~\ref{prop:l0}~(\ref{item:diaminfipropl0})), hence
$x\in\JP=\partial\KP.$

Let $x$ be in the Julia set of $P$ and let $n\in\N.$ From the definition of $L_0$ and Proposition~\ref{prop:l0}~(\ref{item:l0precpl0propl0}) we
have that $L_0$ belongs to $]P^n(x),P^n(L_0)[.$ Hence, the intersection $[x,L_0]\cap\cL_n$ contains exactly one element, denoted by $L_n(x).$

Suppose that $\lim L_n(x)=y\ne x.$ It follows that $x\prec y,$ because $x\lec L_n(x)$ for all $n\ge0.$ From the above we conclude that $y$ is a Julia 
point. Following Lemma~\ref{lemma:segmentfatou} we have that $x\in A(y)\sub\FP,$ which is impossible. Hence we have $\lim L_n(x)=x.$
\end{proof}

\begin{definition}
We will refer to the sequence $(L_n(x))_{n\ge0}$ constructed in the proof of Proposition~\ref{prop:juliasequence} as \emph{the dynamical sequence of
$x$.}
\end{definition}

The dynamical sequences of $x$ and $P(x)$ are related according the following identity
\[L_n(P(x))=P(L_{n+1}(x)),\]
for all $n\ge1.$

As an immediate consequence of Proposition~\ref{prop:juliasequence} we have the following Corollary.

\begin{corollary}
Let $P\in K[z]$ be a polynomial of degree $\ge2.$ Denote by $\hD(P)$ the set of dynamical sequences of $P$ endowed with the topology induced by the
following distance
\[d((L_n),(L_n'))=\frac{1}{m}\]
where $m=\min\{j\ge0\mid L_j\ne L_j'\}$ and $d((L_n),(L_n))=0.$

Let $\what P\colon\hD(P)\to\hD(P)$ be the map defined by $\what P((L_n))=(P(L_{n+1})).$ Then $P\colon\JP\to\JP$ is topologically conjugate to
$\what P\colon\hD(P)\to\hD(P).$ The topological conjugacy is given by $\sigma\colon\JP\to\hD(P)$ where $\sigma(x)=(L_n(x)).$
\end{corollary}

To distinguish whether a Julia point $x$ is classical (i.e $x\in K$) or not (i.e $x\in\HK$) we consider the hyperbolic distance between the
level 0 dynamical point $L_0$ and $x.$ In view of Lemma~\ref{lemma:distlin} we have that
\begin{align*}
\dH(L_n(x),L_0)&=\log(\diam(L_0))-\log(\diam(L_n(x)))\\
&=\sum_{j=0}^{n-1}\dH(L_j(x),L_{j+1}(x)).
\end{align*}
Hence $x$ is a classical point if and only if the sum of the right hand side of the expression above is divergent. The
convergence of the sum, does not allow us to decide whether the point $x$ is of type II, III or IV it only says that
$\log(\diam(x))$ is a positive rational or irrational number.

The following corollaries are applications of Proposition~\ref{prop:finitecritical} to dynamical sequences.

\begin{corollary}
Let $P\in K[z]$ be a tame polynomial of degree $\ge2.$ Then there exists $\cM(P)\in\N,$ only depending on $P,$ such that
\begin{enumerate}
\item if $n\ge\cM(P)$ and $L_n$ is a level $n$ point which is critical, then $L_n=L_n(c)$ for some $c\in\Crit(P)\cap\JP.$

\item $\degl{L_n(c)}=\degl{L_{\cM(P)}(c)}$ for all $n\ge\cM(P)$ and all $c\in\Crit(P)\cap\JP.$
\end{enumerate}
\end{corollary}
\begin{proof}
Note that $w\in\CritI(P)$ belongs to $\KP$ if and only if $[w,\infty[\,\cap\cL_n\ne\vacio$ for all $n\ge1.$ It follows that there exists a smallest
integer
$M_1$ such that if $w\in\CritI(P)$ and $[w,\infty[\,\cap\cL_{M_1}\ne\vacio,$ then $w\in\KP.$

From the definition of the local degree, we have that for each $c\in\Crit(P)\cap\JP$ there exist a smallest integer $M_c$ such that
$\degl c=\degl{L_{M_c}(c)}.$ Consider $M_2=\max\{M_c\mid c\in\Crit(P)\cap\JP\},$ we can consider $\max$ instead $\sup$ by Proposition~\ref{prop:finitecritical}.
Then $\cM(P)=\max\{M_1,M_2\}$ only depends on $P$ and is the smallest integer satisfying (i) and (ii).
\end{proof}

To state and prove the following corollary we need two definition.

\begin{definition}\label{def:o+}
Let $P\in K[z]$ be a nonsimple polynomial of degree $\ge2$ and consider $x\in\JP.$ The \emph{forward orbit of $x$} is the set
\[\cO^+(x)\defi\{P^j(x)\mid j\in\N\}.\]
\end{definition}

\begin{definition}\label{def:wlimite}
Let $P\in K[z]$ be a nonsimple polynomial of degree $\ge2$ and consider $x\in\JP.$ The \emph{$\w$-limit of $x$} is the set
\[\w(x)\defi\Big\{y\in\AK\mid\textnormal{ there exists $(n_j)\sub\N,n_j<n_{j+1}$ and $\lim_{j\to+\infty}P^{n_j}(x)=y$}\Big\}\]
\end{definition}

\begin{corollary}\label{cor:M}
Let $P\in K[z]$ be a tame polynomial of degree $\ge2$ and consider $x$ in the Julia set of $P.$ Then there exists an integer
$N\ge\cM(P),$ depending on $x,$ such that
\begin{enumerate}
\item $\degl x=\degl{L_n(x)}$ for all $n\ge N.$

\item if $n\ge N,\,j\ge1$ and $L_n(P^j(x))$ is critical, then $L_n(P^j(x))=L_n(c)$ for some $c\in\barra{\cO^+(x)}\cap\Crit(P).$
\end{enumerate}
\end{corollary}
\begin{proof}
From the definition of local degree, there exists $N_1\ge\cM(P)$ such that
\[\degl x=\degl{L_{N_1}(x)}.\] 

A critical element $c$ belongs to $\w(x)\setminus\cO^+(x)$ if and only if there exists a increasing sequence $(n_j)$ of integers such 
that the dynamical sequence of $P^{n_j}(x)$ coincides with the dynamical sequence of $c$ at least up to the level $j.$ If 
$c\not\in\w(x)\setminus\cO^+(x)$ there exists an integer $N_c$ such that $L_{N_c}(P^j(x))\ne L_{N_c}(c)$ for all $j\ge1.$ By 
Proposition~\ref{prop:finitecritical} we can consider
\[N_2=\max\{N_c\mid c\in (\Crit(P)\cap\JP)\setminus\w(x)\}.\]

It follows that every $N\ge\max\{N_1,N_2\}$ satisfies (i) and (ii).
\end{proof}

\begin{definition}\label{def:N}
Given $x\in\JP$ we define \emph{the good starting level of $x,$} denoted by $\cN(x),$ as the smallest integer satisfying the two
properties in Corollary~\ref{cor:M}.
\end{definition}

In general, for $x\in\JP$ we want to estimate the distance $\dH(x,L_0).$ But in the practice we estimate the distance $\dH(x,L_{\cN}(x)),$ where
$\cN$ is the good starting level of $x,$ since it is easier to control.

\begin{proposition}\label{prop:hlrec}
Let $P\in K[z]$ be a tame polynomial of degree $\ge2$ and consider $x$ in the Julia set of $P.$ If $x$ belongs to $\HK,$ then the $\w$-limit of $x$
contains at least one critical point of $P$
\end{proposition}
\begin{proof}
Consider $x\in\JP.$ Suppose that $\w(x)\cap\Crit(P)$ is empty, it is enough to show that $x\in K.$ Passing to an iterate if necessary we can suppose
that $x$ has no critical iterates. Since $\w(x)\cap\Crit(P)=\vacio$ and $x$ has no critical iterates, we have that $L_n(P^j(x))$ is noncritical for all $n\ge\cN$ and all $j\ge1,$ where
$\cN=\cN(x)$ is the good starting level of $x.$ Equivalently,
\[L_{n-j}(P^j(x))=P^j(L_n(x))\]
is noncritical provided that $n-j\ge\cN$ (see Definition~\ref{def:N}).

Since, the dynamical level sets $\cL_n$ are finite, for all $n\ge1,$ there exist finitely many intervals of the form $[L_{\cN+1}(y),L_{\cN}(y)]$ with
$y\in\JP.$
Hence, there exist a point $y_0$ in $\JP$ and a strictly increasing sequence $(n_j)$ of dynamical levels larger than $\cN$ such that
\[P^{n_j-\cN}([L_{n_j+1}(x),L_{n_j}(x)])=[L_{\cN+1}(y_0),L_{\cN}(y_0)].\]
The levels $n_j$ are larger than $\cN,$ then by Lemma~\ref{lemma:distgrado} we have that
\[\deg_{L_{n_j}(x)}(P^{n_j-\cN})=\deg_{L_{n_j}(x)}(P)=\degl x.\]
It follows that
\begin{align*}
\dH(x,L_\cN(x))&=\sum_{j=\cN}^\infty\dH(L_{j+1}(x),L_{j}(x))\\
&\ge\sum_{j=0}^\infty\dH(L_{n_j+1}(x),L_{n_j}(x))\\
&=[\degl x]^{-1}\sum_{j=0}^\infty\dH(L_{\cN+1}(y_0),L_{\cN}(y_0))\\
&=+\infty.
\end{align*}
Thus, $x$ belongs to $ K.$
\end{proof}

\section{The Geometric Sequence}\label{sec:geosec}

Let $P\in K[z]$ be a tame polynomial of degree $d\ge2,$ consider $x\in\JP$ and let $L_n=L_n(x)$ for all $n\ge1.$

In order to compute the hyperbolic distance between $x$ and $L_0,$ we want to estimate the distance between two consecutive levels of
the dynamical sequence of $x.$ However, Lemma~\ref{lemma:distgrado} not always applies to relate $\dH(L_{n+1},L_n)$ to the
distance $\dH(P(L_{n+1}),P(L_n))$ because the local degree of $P$ is not necessarily constant in the segment $]L_{n+1},L_n[.$ 

In view of this, to have a better control of the distance between consecutive dynamical levels, we need a finer subdivision of the 
segment joining $x$ to the level 0 dynamical point. This subdivision, that will be the called \emph{the geometric sequence of $x$}, is motivated by
the following propositions about \emph{branch points} of $\CJP.$ First we need a formal definition. 

\begin{definition}
We say that a type II point $x\in\CJP$ is a \emph{branch point of $\CJP$} if $\CJP$ intersects at least two directions in $T_x.$
\end{definition}

\begin{proposition}
Let $P\in K[z]$ be a tame polynomial of degree $\ge2.$ Then $L_0$ is a branch point of $\CJP.$
\end{proposition}
\begin{proof}
From Proposition~\ref{prop:l0}~(\ref{item:l0precpl0propl0}) it follows that the elements of $P^{-1}(L_0)$ are strictly smaller than $L_0.$ Following
Proposition~\ref{prop:l0}~(\ref{item:l0branchpropl0}) there exist $x_1,x_2\in P^{-1}(L_0)$ such that $x_1,x_2$ belong to different directions in $T_{L_0}.$
Hence, $L_0=x_1\vee x_2$ is a branch point of $\CJP.$ 
\end{proof}

\begin{proposition}\label{prop:noramas}
Let $P\in K[z]$ be a tame polynomial of degree $\ge2$ and consider $x_B$ a branch point of $\CJP.$ Then there exists a direction $\cD\in T_{x_B}$ with
$\cD\cap\JP=\vacio$ and a critical point $w\not\in\KPI$ such that $P^n(w)\in\cD$ for some $n\ge0.$
\end{proposition}
\begin{proof}
In view of Proposition~\ref{prop:juliasequence} there exits $n$ such that $L_n\prec x_B\lec L_{n-1}$ for some dynamical points $L_n$ and
$L_{n-1},$ it follows that $P^n(x_B)$ belongs to the segment $]L_0,P(L_0)].$

Since the dynamical points of a given dynamical level are finitely many, we have that there exist only finitely many directions in $T_{x_B}$ having nonempty
intersection with $\JP.$ Let $\cD_1,\dots,\cD_m$ be such directions. 

We have that $P^n(\cD_j)=\cD_0$ for all $1\le j\le m,$ where $\cD_0$ is the direction at $y$ containing $L_0.$ From the Riemann-Hurwitz
formula (see Remark~\ref{rem:RH}), we obtain that
\begin{align*}
\deg_{x_B}(P^n)&=1+\sum_{z\in B\cap\CritI(P^n)}(\deg_z(P^n)-1)\\
&=\deg(T_{x_B}P^n)\\
&=\deg_{\cD_1}(P^n)+\cdots+\deg_{\cD_m}(P^n)\\
&=m+\sum_{z\in(\cD_1\sqcup\cdots\sqcup\cD_m)\cap\CritI(P^n)}(\deg_z(P^n)-1).
\end{align*}
Since $2\le m,$ it follows that
\[\sum_{z\in(\cD_1\sqcup\cdots\sqcup\cD_m)\cap\CritI(P^n)}(\deg_z(P^n)-1)<\sum_{z\in B\cap\CritI(P^n)}(\deg_z(P^n)-1).\]
That is, there exists a critical point of $P^n$ in a direction $\cD\in T_{x_B}\setminus(\cD_1\sqcup\cdots\sqcup\cD_m).$
Hence, the proposition follows.
\end{proof}

Let $\{w_1,\dots,w_{q'}\}$ be the set of classical critical points of $P$ which are not contained in $\KPI.$ Note that $1\le q'\le d-1.$ 
For each $1\le j\le q',$ let $x_{w_j}$ be defined by
\[(\CJP\cup[L_0,\infty[)\cap[w_j,\infty[\,=[x_{w_j},\infty[.\]

By Proposition~\ref{prop:noramas} all the branch points of $\CJP$ are contained in the grand orbit of the set $\{x_{w_1},\dots,x_{w_{q'}}\}.$

\begin{definition}\label{def:geo}
\emph{The geometric sequence of $x,$} denoted by $(G_n(x))_{n\ge0},$ is the decreasing sequence enumerating the elements in
\[[x,L_0]\cap\GO(\{x_{w_1},\dots,x_{w_{q'}}\}),\]
where $\GO$ denotes the grand orbit.
\end{definition}

\begin{remark}\label{rem:despues}
Consider $v_{q-1},\dots,v_0\in]L_0,P(L_0)]$ satisfying
\[\GO(\{x_{w_1},\dots,x_{w_{q'}}\})\,\cap\,]L_0,P(L_0)]=\{v_{q-1},v_{q-2},\dots,v_0\},\]
and $v_{q-1}\prec v_{q-2}\prec\cdots\prec v_0=P(L_0).$ Note that $1\le q\le d$ and that every element of
\[\GO(\{x_{w_1},\dots,x_{w_{q'}}\})\cap\CJP\]
is eventually mapped to $v_j,$ for some $j.$

Consider $n\ge0$ and put $\displaystyle{n=q\cdot\pen{n/q}+j}$ for some $0\le j\le q-1,$ where $\pen\cdot$ denotes the floor function, that is $\pen\cdot\colon\R\to\Z$
is defined by
\[\pen a\defi\sup\{n\in\Z\mid n\le a\}.\] 

It follows that
\[P^{\pen{n/q}+1}(G_n(x))=v_j.\]
Moreover
\[P^{\pen{n/q}+1}([G_{n+1}(x),G_{n}(x)])=I_j,\]
where $I_j=[v_{j+1},v_{j}]$ for $0\le j\le q-2$ and $I_{q-1}=[L_0,v_{q-1}].$
\end{remark}

\begin{definition}
We say that $v_{q-1},\dots,v_0$ are the \emph{generators of the geometric sequences of $P.$}
\end{definition}

The next lemma states the basic properties of the geometric sequence. The proof is straightforward and we omit it.

\begin{proposition}
Let $P\in K[z]$ be a tame polynomial of degree $\ge2$ with $q$ generators for its geometric sequences and consider $x\in\JP.$ Then, the following statements holds 
\begin{enumerate}
\item $L_n(x)=G_{nq}(x)$ for all $n\ge 0.$

\item $\displaystyle{\lim_{n\to+\infty}G_n(x)=x}.$ 

\item The geometric sequences of $x$ and $P(x)$ are related according
\[P(G_n(x))=G_{n-q}(P(x))\]
for all $n\ge q.$
\end{enumerate}
\end{proposition}

Note that
\[P(]G_{n+1}(x),G_{n}(x)[)=]G_{n-q+1}(P(x)),G_{n-q}(P(x))[\]
for all $n\ge q+1.$ Moreover, for all $n\ge0$ it follows that
\[P^{-1}(]G_{n+1}(x),G_{n}(x)[)=\bigcup_{y\in P^{-1}(x)}]G_{n+q+1}(y),G_{n+q}(y)[\]

The main result in this section is the following proposition. It will allow us to use Lemma~\ref{lemma:distgrado} to relate $\dH(G_{n+1},G_n)$ with $\dH(P(G_{n+1}),P(G_n)).$

\begin{proposition}\label{prop:enecalig}
Let $P\in K[z]$ be a tame polynomial of degree $\ge2$ and consider $x\in\JP.$ Let $(G_n)$ be the geometric sequence of $x.$ Then for all $n\ge0$
the local degree of $P$ is constant in the segment $[G_{n+1},G_{n}[.$ In fact
\[\degl y=\degl{G_{n+1}}\]
for all $y\in[G_{n+1},G_{n}].$
\end{proposition}
\begin{proof}
Let $n\ge0$ and define
\[\mathcal I_n\big\{y\in\,]G_{n+1},G_{n}]\mid\degl w=\degl{G_{n+1}}\textnormal{ for all } w\in[G_{n+1},y[\big\}.\]
To show that $\mathcal I_n\ne\vacio,$ consider $a\in K$ and $r>0$ such that $B_r^+(a)$ is the ball associated to $G_{n+1}.$ Consider
\[R=\min\{\abs{a-w}\mid w\in\CritI(P)\setminus B_r^+(a)\}.\]
In view of the Riemann-Hurwitz formula, see Remark~\ref{rem:RH}, the local degree of $P$ at the point associated to the ball $B_{r+\e}^+(a)$ coincides
with $\degl{G_{n+1}}$ for all $0\le\e< R.$ In particular, we obtain that $\mathcal I_n\ne\vacio.$ Moreover, the previous argument shows that
largest element contained in $\mathcal I_n$ is a type II point.

Denote by $y_B=\max\,\mathcal I_n.$ We will show that $y_B=G_{n}.$

We proceed by contradiction. Suppose that $y_B\prec G_{n}.$ It follows that the degree of the map $P$ in $\barra B$ is larger than the degree
of $P$ in the direction $\cD$ at $y_B$ that contains $G_{n+1}.$ In particular, there exists a critical point $w$ in $\barra B\setminus\cD.$

We have two cases: 

If $w\not\in\KPI,$ then $y_B$ belongs to $\{x_{w_1},\dots,x_{w_{q'}}\}$ (see Definition~\ref{def:geo}). In particular, $G_{n}$ and $G_{n+1}$ are
not consecutive elements of the sequence $(G_n).$

If $w\in\KPI$ we have that $y_B$ is a branch point of $\JP,$ which is a contradiction because $]G_{n+1},G_{n}[$ is branch point free, see
Proposition~\ref{prop:noramas}.

It follows that $\max\mathcal I_n=G_{n}.$ Therefore $\degl y=\degl{G_{n+1}}$ for all $y$ in $[G_{n+1},G_{n}[.$
\end{proof}

\section{Algebraic Degree on the Affine Line}\label{sec:algdeg}

From now, we consider that $K$ is an algebraically closed field which is complete with respect to a non-Archimedean absolute value. We will assume that
there exits a complete field $k\sub K$ so that $\abs\cdot$ restricted to $k$ is a discrete absolute value and such that the elements of $K$ which are algebraic over $k$ are dense in $K,$ that is
\[\{z\in K\mid [k(z):k]<+\infty\}\]
is a dense subset of $K.$ Moreover, changing $\abs\cdot$ for $\abs\cdot^\lambda$ for some $\lambda>1$ we can always suppose that
\[\log(\abs{k^\times})=\Z\quad\textnormal{and}\quad\log(\abs{K^\times})=\Q.\]
Since $K$ is algebraically closed an complete, it is not difficult to see that $K$ coincides with the completion of an algebraic closure of $k.$

The algebraic degree of an element $z\in K$ over $k$ is the number $[k(z):k],$ that is, the degree of the smallest extension of $k$ containing $z.$ To extend this notion to the Berkovich line, note that every nonsingular element in $\HK$ belongs to the convex hull of $\calg k\cap K,$ where $\calg k$ is an algebraic closure of $k.$ It follows that, if $x$ is a type II or III point
in $\HK,$ then
\[\{[F:k]\mid x\in\conv(F), F\textnormal{ is a finite extension of $k$}\}\ne\vacio,\]
where $\conv(F)$ denotes the convex hull of $F,$ see subsection~\ref{subsec:berk}.

\begin{definition}
For $z\in K$ we define $\delta(z)\defi [k(z):k].$ For $x\in\AK$ nonsingular, we define the \emph{algebraic degree of $x$ over $k$} as
\[\delta(x)=\min\{[F:k]\mid x\in\conv(F), F\textnormal{ is a finite extension of $k$}\}.\]
If $x\in\AK$ singular, we define the algebraic degree of $x$ as $+\infty.$
\end{definition}

\begin{definition}
We say that a nonempty subset $X$ of $\AK$ is an \emph{algebraic set over $k$} if $\delta(x)<+\infty$ for all $x\in X.$
\end{definition}

Note that for $x_B\in\HK$ nonsingular and $F$ a finite extension of $k,$ we have that $x_B$ belongs to $\conv(F)$ if and only
if $B\cap F\ne\vacio.$ In fact, if $x_B$ belongs to $\conv(F),$ by definition there exist two points $z_0,z_1\in F$ such that $x_B$ belongs to the
segment $]z_0,z_1[.$ In particular, $z_0\prec x_B$ or $z_1\prec x_B,$ that is $z_0\in B\cap F$ or $z_1\in B\cap F.$ Conversely, if 
$B=B_r^+(z_0)$ with $z_0\in F$ then $z_0\lec x_B.$ It follows that $x_B$ belongs to $\conv(F).$ Hence, the algebraic degree of a nonsingular element
$x_B\in\HK$ is
\[\delta(x_B)=\delta(B)\defi\min\{[k(a):k]\mid a\in B\}.\]

The next lemma states the basic properties of the algebraic degree. The proof is straightforward
and we omit it.

\begin{lemma}\label{lemma:justi}
Let $P$ be a polynomial with coefficients in $k$ of degree $\ge2.$ Then the following statements hold:
\begin{enumerate}
\item $\delta(P(x))\le\delta(x)$ for all $x\in\AK.$\label{item:gradobajalemmajusti}

\item If $x\lec y,$ then $\delta(y)\le\delta(x).$\label{item:leclemmajusti}

\item $\displaystyle{\delta(x)=\lim_{n\to+\infty}\delta(x_n)}$ for each decreasing sequence $(x_n)$ contained in $\HK$ such that $x_n\to x.$\label{item:gradoseclemmajusti}
\end{enumerate}
\end{lemma}

In the case of the Julia set, the fact of being an algebraic set over $k$ is a local property:

\begin{proposition}
Let $P\in k[z]$ be a nonsimple polynomial of degree $\ge2.$ The Julia set of $P$ is algebraic over $k$ if and only there exists a point $x\in\JP$ and a
neighborhood $V$ of $x$ such that $\JP\cap V$ is an algebraic set over $k.$
\end{proposition}
\begin{proof}
If $\JP$ is algebraic over $k,$ then taking $V$ as an open ball of the Berkovich line which contains $\JP,$ we have that $V$ is a neighborhood of
all $x\in\JP$ and $\JP=V\cap\JP$ is an algebraic set over $k.$

Conversely, if there exists a element $x\in\JP$ and a neighborhood $V$ containing $x$ we have that
\[\JP=\bigcup_{n\ge 0}P^n(\JP\cap V)\sub\bigcup_{n\ge 0}P^n(V).\]
From Lemma~\ref{lemma:justi}~(\ref{item:gradobajalemmajusti}) it follows that $\delta(P^{n+1}(x))\le\delta(P^n(x))$ for all $x\in\JP\cap V$
and all $n\ge0.$ That is, the Julia set of $P$ is an algebraic set over $k.$
\end{proof}

\subsection{Equilibrium measure and algebraic degree}

Let $p\in k[z]$ be a nonsimple polynomial of degree $\ge 2.$ In this subsection, we will show that if there exists a element $y\in\JP$ with $\delta(y)=+\infty$ then $\delta(x)=+\infty$ for $\rho_P$-almost all $x\in\JP.$ 

\begin{lemma}\label{lemma:Delta}
Let $P\in k[z]$ be a nonsimple polynomial of degree $\ge 2$ with $\JP$ algebraic over $k.$ Let $\cA_n=\{x\in\JP\mid\delta(x)\ge n\}.$ If $\cA_n\ne\vacio$ then $\rho_P(\cA_n)=1$
where $\rho_P$ is the equilibrium measure of $P.$
\end{lemma}
\begin{proof}
If there exists $x\in\JP$ with $\delta(x)\ge n$ then there exists $m_n\in\N$
such that $\delta(L_{m_n}(x))\ge n,$ see Lemma~\ref{lemma:justi}~(\ref{item:gradoseclemmajusti}). Let $B_n$ be the closed ball 
associated to $L_{m_n}.$ We have that $\rho_P(\barra{B_n})>0,$ since $\barra{B_n}$ contains an open set that contains a Julia point, then by the
ergodicity of $\rho_P$ it follows that $\rho_P$-almost all $x\in\AK$ have infinitely many iterates in $\barra{B_n}.$ From 
Lemma~\ref{lemma:justi}~(\ref{item:gradobajalemmajusti}) and (\ref{item:leclemmajusti}) we conclude that $\rho_P(\cA_n)=1.$
\end{proof}

\begin{proposition}\label{prop:Delta}
Let $P\in k[z]$ be a nonsimple polynomial of degree $\ge2$ with $\JP$ algebraic over $k.$ Then the set $\{\delta(x)\mid x\in\JP\}$ is bounded. Moreover,
\[\rho_P(x\in\JP\mid\delta(x)=\Delta)=1\]
where $\Delta=\max\{\delta(x)\mid x\in\JP\}$ and $\rho_P$ is the equilibrium measure of $P.$
\end{proposition}
\begin{proof}
For each $n\in\N$ let $\cA_n$ as in Lemma~\ref{lemma:Delta}.

Suppose that $\{\delta(x)\mid x\in\JP\}$ is an unbounded set. Then there exists a sequence $(x_n)$ of Julia points such that $\delta(x_n)\ge n$
and therefore, using Lemma~\ref{lemma:Delta}, it follows that $\rho_P(\cA_n)=1$ for all $n\in\N.$ Since  $\cA_{n+1}\sub\cA_{n}$ and
\[\displaystyle{\rho_P(\cap\cA_n)=\lim_{n\to+\infty}\rho_P(\cA_n)=1},\]
we have that the intersection $\cap\cA_n$ is nonempty. It follows that there exists $y\in\JP$ with $\delta(y)=+\infty$, which is a contradicts our 
assumption that $\JP$ is an algebraic set over $k.$ Hence, the set $\{\delta(x)\mid x\in\JP\}$ is bounded.

Consider $\Delta=\max\{\delta(x)\mid x\in\JP\}.$ By Lemma~\ref{lemma:Delta} we have that
\[\{x\in\JP\mid\delta(x)=\Delta\}=\cA_{\Delta}\aev\]
Then $\rho_P(\{x\in\JP\mid\delta(x)=\Delta\})=1.$
\end{proof}

The following corollary will be useful in the proofs of Theorem~\ref{theo:juliaalg} and of Theorem~\ref{theo:juliaalgiff}.

\begin{corollary}\label{cor:unotodo}
Let $P\in k[z]$ be a polynomial of degree $\ge2$ and denote by $\rho_P$ the equilibrium measure of $P.$ Then the following statements are
equivalent:
\begin{enumerate}
\item There exists $y\in\JP$ with $\delta(y)=+\infty.$

\item For each $n\ge1$ there exists $y_n\in\JP$ such that $\delta(y_n)\ge n.$

\item $\rho_P(x\in\JP\mid\delta(x)=+\infty)=1.$
\end{enumerate}
\end{corollary}

\section{Polynomial with Algebraic Coefficients}\label{sec:polalg}

In this section we will fix a polynomial $P$ with algebraic coefficients over $k$ and we will study how the algebraic degree behaves along the geometric sequence of $x$ in $\JP$. To do this, we need the following \emph{dynamical version} of the well known Krasner's Lemma which is adapted for our applications. See Corollary 3 in chapter seven of~\cite{local_cassels} for the standard version of the lemma.

\begin{lemma}[Krasner's Lemma]\label{lemma:krasner}
Let $k$ be a field of characteristic $0$ and complete with respect to a non-Archimedean absolute value. Let $\calg k$ be an algebraic closure of the
field $k.$
Consider $P\in k[z]$ and let $\alpha\in\calg k$ such that $P (\alpha) = 0$. If $B\sub\calg k$ is a ball containing $\alpha$ such that $\degl B=1$,
then
$k(\alpha)\sub k(\beta)$, for all $\beta\in B.$
\end{lemma}

\subsection{Algebraic degree along a geometric sequence}

The first statement of Lemma~\ref{lemma:justi} shows that the algebraic degree behaves nicely under the action of polynomials with algebraic coefficients. 

Through this subsection, let $P$ be a tame polynomial with algebraic coefficients over $k.$ Passing to a finite extension if necessary, we can suppose that the coefficients and the critical points belongs to $k.$
\vspace{.4cm}

To state the proposition that allow us to estimate the algebraic degree along a geometric sequence we need two definitions.

\begin{definition}
Consider $x$ in the Julia set of $P,$ and let $(G_n)$ be its geometric sequence. For each $n\ge0$
\emph{the injectivity time of $G_n$} is the largest integer $0\le t_n$ such that $\deg_{G_n}(P^{t_n})=1.$
The \emph{critical pullbacks of $P$ around $G_n$} are the elements in $D_n\cap P^{-t_n}(\CritI(P)),$ where $x_{D_n}=G_n.$
\end{definition}

\begin{definition}\label{def:sn}
Let $P\in k[z]$ be a nonsimple polynomial of degree $\ge2$ and consider $x\in\JP.$ Let $t_{-1}=-1$ and $s_{-1}=1.$ For each $n\ge0$ define $s_n(x)$ as the index of $\abs{k^\times}$
in the group generated by $\abs{k^\times}$ and $\diam(G_n(x)),$ that is
\[s_n(x)\defi\Big[\abs{k^\times}(\diam(G_n(x))):\abs{k^\times}\Big]\]
\end{definition}

\begin{proposition}\label{prop:condicion}
Let $P\in k[z]$ be a tame polynomial of degree $\ge2.$ Consider $x$ in the Julia set of $P$ and denote by
$(G_n)$ the geometric sequence of $x.$ For each $n\ge-1$ let $s_n=s_n(x).$ Then
\begin{equation}\label{eq:condicion}
\max\{s_j\mid t_j\ne t_{j+1},-1\le j\le n\}\le\delta(G_{n+1}).
\end{equation}
\end{proposition}

The reader may find the proof of this proposition is at the end of this subsection.

Note that if $P$ has $q$ generators for its geometric sequences then the injectivity time
$t_n$ is smaller or equal than $\pen{n/q}+1.$ In the case that $G_n$ is a critical element we have that $t_n=0$ and that the critical pullbacks around $G_n$
are the critical points of $P$ contained in $D_n.$

If $t_n>0$ then the injectivity time of $G_n$ is the smallest integer such that
\[P^{t_n}(G_n)\in\Crit(P)\quad\textnormal{and}\quad P^{t_n-1}(G_n)\not\in\Crit(P).\]

\begin{lemma}\label{lemma:degree}
Let $P\in k[z]$ be a tame polynomial of degree $\ge2.$ Consider $x\in\JP$ and denote by $(G_n)$ its geometric
sequence.
Consider critical pullbacks $u_n,u_{n+1}$ around $G_n$ and $G_{n+1}$ respectively. Then, for all $n\ge0$ the followings statements hold:
\begin{enumerate}
\item $\delta(G_n)=\delta(u_n)=\delta(u_n')$ for all critical pullbacks $u_n'$ around $G_n.$\label{item:gradopreimlemmadegree}

\item $k(u_n)$ is contained in $k(u_{n+1}).$\label{item:extefinlemmadegree}

\item $u_{n}\not\in\cD_n$ and $\delta(G_{n+1})=\delta(\cD_n\cap K),$ where $\cD_n$ is the direction at $G_n$ that
contains $G_{n+1}.$\label{item:distvnlemmadegree}
\end{enumerate}
\end{lemma}
\begin{proof}
Note that if $G_n$ is a critical element we have that $u_n$ is a critical point of $P.$ Therefore (\ref{item:gradopreimlemmadegree}) and
(\ref{item:extefinlemmadegree}) follow directly, in this case.

If $G_n$ is noncritical, then $t_n>1$ and
\[\deg_{D_n}\big[P^{t_n}(z)-P^{t_n}(u_n)\big]=1.\]
By Krasner's lemma (see Lemma~\ref{lemma:krasner}), it follows that $k(u_n)\sub k(v)$ for all $v$ in $D_n\cap\fpui.$

Taking $v\in D_n$ such that $\delta(v)=\delta(G_n)$ we have that $\delta(v)\le\delta(u_n)\le\delta(v),$ that is $\delta(u_n)=\delta(G_n).$ If $u_n'$
is
another critical pullback around $G_n$ then interchanging $u_n$ with $u_n'$ we obtain the second equality in (\ref{item:gradopreimlemmadegree}).

Taking $v=u_{n+1}$ we obtain (\ref{item:extefinlemmadegree}).

In order to prove (\ref{item:distvnlemmadegree}), recall that according to Proposition~\ref{prop:enecalig} we have that
\[\deg_{G_{n+1}}(P^{j})=\deg_{y}(P^{j})\]
for all $y\in[G_{n+1},G_n[$ and all $j\ge1.$ In particular $t_{n+1}$ is the smallest integer such that $P^{t_{n+1}}(\cD_n)\cap\Crit(P)\ne\vacio$ and 
$P^{t_{n+1}-1}(\cD_n)\cap\Crit(P)=\vacio.$
Then, as in the proof of the first part of the lemma, we have that
\[\deg_{\cD_n\cap K}\big[P^{t_{n+1}}(z)-P^{t_{n+1}}(u_{n+1})\big]=1.\]
Therefore $\delta(G_{n+1})=\delta(\cD_n\cap K).$
\end{proof}

\begin{corollary}\label{cor:tndifferent}
Let $P\in k[z]$ be a tame polynomial of degree $\ge2.$ Let $x$ be a noncritical Julia point and denote
by $(G_n)$ its geometric sequence.
\begin{enumerate}
\item If $t_{n}=t_{n+1}$ then $\delta(G_n)=\delta(G_{n+1}).$

\item If $\delta(G_n)<\delta(G_{n+1})$ then $t_n<t_{n+1}.$
\end{enumerate}
\end{corollary}

Now we can give the proof of Proposition~\ref{prop:condicion}. This proposition is a key ingredient in order to prove our main results.

\begin{proof}[Proof of Proposition~\ref{prop:condicion}]
We proceed by induction in $n.$ For $n=-1$ we have that $-1=t_{-1}\ne t_0=0,$ since $\degl{L_0}=\deg(P),$ therefore
\[\max\{s_j\mid t_j\ne t_{j+1},-1\le j\le-1\}=\max\{s_{-1}\}=1\le\delta(G_0)=1.\]
That is, the first step of the induction is valid.

Suppose that (\ref{eq:condicion}) is valid for $n-1.$ If $n$ is such that $t_{n}=t_{n+1}$ we have that,
$\{s_j\mid t_j\ne t_{j+1},-1\le j\le n-1\}=\{s_j\mid t_j\ne t_{j+1},-1\le j\le n\}.$ It follows 
that $\max\{s_j\mid t_j\ne t_{j+1},-1\le j\le n\}\le\delta(G_n)\le\delta(G_{n+1}).$

Suppose that $n$ is such that $t_{n}<t_{n+1}.$ In this case, the corresponding
critical pullbacks $u_n$ and $u_{n+1}$ are necessarily different. By Lemma~\ref{lemma:degree}~(\ref{item:extefinlemmadegree}), we have that
$k(u_{n+1})$
is a finite extension of $k(u_n),$ and therefore $\delta(u_n)$ divides $\delta(u_{n+1}).$

Following Lemma~\ref{lemma:degree}~(\ref{item:distvnlemmadegree}) we have that the distance between $u_n$ and $u_{n+1}$ coincides with the diameter
of $G_n.$ Therefore
\begin{align*}
s_n&=\Big[\abs{k^\times}(\diam(G_n)):\abs{k^\times}\Big]\\
&=\Big[\abs{k^\times}(\abs{u_n-u_{n+1}}):\abs{k^\times}\Big]\\
&\le[k(u_n-u_{n+1}): k]\\
&\le\max\{\delta(u_n),\delta(u_{n+1})\}.
\end{align*}
Hence, we obtain that
\[\max\{s_n,\delta(u_n)\}\le\delta(u_{n+1}).\]
Applying the inductive hypothesis it follows that
\begin{align*}
\delta(u_{n+1})&\ge\max\{s_n,\delta(u_n)\}\\
&\ge\max\{s_n,\max\{s_j\mid t_j\ne t_{j+1},-1\le j\le n-1\}\}\\
&=\max\{s_j\mid t_j\ne t_{j+1},-1\le j\le n\}.
\end{align*}
Thus, we have proven Proposition~\ref{prop:condicion}.
\end{proof}

\subsection{No wandering components for polynomials with algebraic coefficients}

The following proposition is the key to prove Corollary~\ref{cor:nowandering}. In fact, combining the proposition below with
Proposition~\ref{prop:fatoujulia} we obtain Corollary~\ref{cor:nowandering} in the case of polynomials with algebraic coefficients over $k.$

\begin{proposition}\label{prop:algnowand}
Let $P\in k[z]$ be a tame polynomial of degree $\ge 2.$ If $x$ is a nonpreperiodic algebraic element in $\JP,$
then $x\in K.$
\end{proposition}
\begin{proof}
Let $x$ be a noncritical and nonpreperiodic algebraic element of the Julia set of $P.$ Since $x\in\JP$ if and only if 
$P^j(x)\in\JP$ for all $j\ge1$ and there are only finitely many critical elements in $\JP,$ we may assume that the forward orbit of $x$ is free of critical elements.
Denote by $(G_n)$ the geometric sequence of $x.$

Since $x$ is not the preimage of a critical element it follows that
\[\{n\in\N\mid t_n\ne t_{n+1}\}\]
contains infinitely many elements. For each $n\ge1$ let $m_n$ be the $n$-th nonnegative integer such that $t_{m_n}\ne t_{m_n+1}.$

In view of Proposition~\ref{prop:condicion}, the set $\{s_{m_j}(x)\mid j\in\N\}$ is bounded by $\delta(x)$. Hence, the denominators of $\log(\diam(G_{m_n}))$ are bounded. 

Let $D\in\N$ be the maximum between the denominators
of $\log(\diam(G_{m_n})),\,n\in\N,$  and $\log(\diam(L_0)).$ It follows that
\[\dH(G_{m_j},G_{m_{j-1}})=\log(\diam(G_{m_{j-1}}))-\log(\diam(G_{m_j}))\ge\frac{1}{D!}\]
for all $j\ge1.$ Hence
\[\dH(G_{m_n},L_0)=\dH(G_{m_0},L_0)+\sum_{j=1}^{n}\dH(G_{m_{j}},G_{m_{j-1}})\ge\dH(G_{m_0},L_0)+\frac{n}{D!}.\]
Thus,
\begin{align*}
d(x,L_0)&=\lim_{n\to+\infty}\dH(G_n,L_0)\\
&=\lim_{n\to+\infty}\dH(G_{m_n},L_0)\\
&=\lim_{n\to+\infty}\Big(\dH(G_{m_0},L_0)+\frac{n}{D!}\Big)\\
&=+\infty.
\end{align*}
Therefore, $x$ belongs to $\JPI.$
\end{proof}

\subsection{Algebraic Julia sets and recurrent critical elements}

Through this subsection let $P\in k[z]$ be a tame polynomial of degree $\ge2$ and $q$ generator for its geometric sequences. Recall
that we are assuming that $k$ contains the critical points of $P.$ 

Consider $x\in\JP$ and let $G_n=G_n(x)$ for all $n\ge1.$ From Proposition~\ref{prop:condicion} it follows that there exists a relation between
$\diam(G_n)$ and $\delta(G_n).$ Moreover. note that we can obtain $\diam(G_n)$ from $\dH(G_n,L_0).$ 
 
From Lemma~\ref{lemma:distgrado} it follows that for each $n\ge0$ the lengths of the
segments $[G_{n+1},G_{n}]$ and $P^{\pen{n/q}+1}([G_{n+1},G_{n}])=I_{j}$ (see Remark~\ref{rem:despues})
are related according the following identity
\begin{equation}\label{eq:eqfea}
\dH(G_{n+1},G_{n})=[\deg_{G_{n+1}}(P^{\pen{n/q}+1})]^{-1}\cdot\abs{I_j},
\end{equation}
where $\abs{I_j}$ denotes the length of the segment $I_j$ for all $0\le j\le q-1.$
Therefore
\[\sum_{n=0}^{+\infty}\dH(G_{n+1},G_{n})=\sum_{n=0}^{+\infty}[\deg_{G_{n+1}}(P^{\pen{n/q}+1})]^{-1}\cdot\abs{I_{j_n}},\]
where, $0\le j_n\le q-1$ is such that $P^{\pen{n/q}+1}([G_{n+1},G_{n}])=I_{j_n}.$

In order to write the previous sum in a more convenient manner we need to introduce the following notation.

\begin{definition}
Consider $x\in\JP.$ For all $n\ge0$ we define, the \emph{dynamical degree of level $n$ around $x$} as
\[d_n(x)\defi\deg_{G_{n+1}(x)}(P^{\pen{n/q}+1}).\]
\end{definition}

\begin{definition}
The \emph{range of dynamical degrees of $x$} is defined by
\[\cD(x)\defi\{d_n(x)\mid n\ge0\}.\]
\end{definition}

Rephrasing equation (\ref{eq:eqfea}) in this notation gives the following lemma.

\begin{lemma}\label{lemma:sumafea}
Let $P\in K[z]$ be a tame polynomial of degree $\ge2$ with $q$ generators for its geometric sequences. Consider $x\in\JP$ and denote by
$(G_n)$ the geometric sequence of $x.$ Then
\[\dH(G_{n+1},L_0)=\sum_{\ell=0}^n\dH(G_{\ell+1},G_{\ell})=\sum_{j=0}^{q-1}\Big(\abs{I_j}\sum_{\substack{\ell\equiv j\bmod{q} \\ 0\le\ell\le n}}[d_\ell^{-1}]\Big)\]
for all $n\ge0,$ where $\abs{I_j}$ is the length of the interval $I_j,\,0\le j\le q-1$
\end{lemma}

The above sum gives us a relation between the sum of the length of the intervals defined by the generators of the geometric
sequences of $P$ and the range of dynamical degrees around $x.$

\subsection{Proof of Theorem~\ref{theo:juliaalg}}

In order to prove Theorem~\ref{theo:juliaalg} and Theorem~\ref{theo:juliaalgiff} we need to establish some relations between the range of dynamical degrees and the 
existence of recurrent critical points.

First we need to recall some notation. Let $P\in K[x]$ be a tame polynomial, see Definition~\ref{def:tame}, and consider $x$ in the Julia set of $P.$ We denote by $\w(x)$ and $\cO^+(x)$ the $\w$-limit and the forward
orbit of $x$ respectively. We denote by $G_n=G_n(x)$ the geometric sequence of $x,$ see Definition~\ref{def:geo}. The fundamental property of the geometric sequence
is the fact that
\[\degl y=\degl{G_n}\]
for all $y\in[G_n,G_{n-1}[.$ 

Recall that, for $a\in\N$ and $X$ a subset of $\N,$ the set $a\cdot X$ is defined by
\[a\cdot X\defi\{a\cdot b\mid b\in X\}.\]

\begin{lemma}\label{lemma:sigma}
Let $P\in K[z]$ be a tame polynomial of degree $d\ge2$ and let $x\in\JP.$ Then the following statements hold:
\begin{enumerate}
\item If $\w(x)\cap\Crit(P)=\vacio,$ then $\cD(x)$ is finite.

\item If $\w(x)\cap\Crit(P)\ne\vacio,$ then the range of dynamical degrees $\cD(x)$ is contained in
\[\{1,2,\dots,d^{\cN(x)+2}\}\cup\Bigg(\degl x\cdot\bigcup_{c\in\cC(x)}\cD(c)\Bigg),\]
where $\cC(x)=\Crit(P)\cap\barra{\cO^+(x)}.$
\label{item:unionlemmasigma}
\end{enumerate}
\end{lemma}
\begin{proof}
Suppose that $P$ has $q$ generators for its geometric sequences and let $\cN=\cN(x)$ be the good starting level of $x$ (see Definition~\ref{def:N}).

Note that if $0\le n\le q\cN,$ then 
\[d_n=\deg_{G_{n+1}}(P^{\pen{n/q}+1})\le \deg_{G_{n+1}}(P^{\pen{q\cN/q}+1})\le d^{\cN+1}.\]

In order to prove the first statement observe that if $\w(x)\cap\Crit(P)=\vacio$ then $x$ has at most $d-2$ critical images.
Consider $n>q\cN$ and let $\ell\ge1$ be the smallest integer such that $(n+1)-\ell q\le q\cN.$ It follows that,
\begin{align*}
d_n&=\deg_{G_{n+1}}(P^{\pen{n/q}+1})\\
&=\deg_{G_{n+1}}(P^\ell)\cdot\deg_{P^\ell(G_{n+1})}(P^{\pen{n/q}+1-\ell})\\
&=\deg_{G_{n+1}}(P^\ell)\cdot\deg_{G_{n-q\ell+1}(P^\ell(x))}(P^{\pen{n/q}+1-\ell})\\
&=\deg_{G_{n+1}}(P^\ell)\cdot\deg_{G_{n-q\ell+1}(P^\ell(x))}(P^{\pen{(n-q\ell)/q}+1})\\
&\le d^{d-1}\cdot d^{\cN+1}\\
&\le d^{d+\cN}.
\end{align*}
Therefore, the range of dynamical degrees of $x$ is contained in the finite set 
\[\{1,2,\dots,d^{d+\cN}\}.\]

To prove (\ref{item:unionlemmasigma}), consider $n>q\cN.$ Then there exists a minimal $1\le j_n\le\pen{n/q}+1$
such that $P^{j_n}(G_{n+1})$ is critical element.

We have two cases:

If $n-qj_n\le q\cN$ we have that
\begin{align*}
d_n&=\deg_{G_{n+1}}(P^{\pen{n/q}+1})\\
&=\deg_{G_{n+1}}(P^{j_n})\cdot\deg_{P^{j_n}(G_{n+1})}(P^{\pen{n/q}+1-j_n})\\
&=\degl{x}\cdot\deg_{G_{n-qj_n+1}}(P^{\pen{(n-qj_n)/q}+1})\\
&\le\degl x\cdot d^{\cN+1}\\
&\le d^{\cN+2}.
\end{align*}
Therefore, $d_n$ belongs to $\{1,2,\dots,d^{\cN+2}\}.$

Suppose that $n-qj_n>q\cN.$ We have that $P^{j_n}(G_{n+1})=G_{n-qj_n+1}(c)$ for some critical element $c\in\Crit(P)\cap\barra{\cO^+(x)}$ (see
Corollary~\ref{cor:M}). It follows that
\begin{align*}
d_n&=\deg_{G_{n+1}}(P^{\pen{n/q}+1})\\
&=\deg_{G_{n+1}}(P^{j_n})\cdot\deg_{P^{j_n}(G_{n+1})}(P^{\pen{n/q}+1-j_n})\\
&=\degl{x}\cdot\deg_{G_{n-qj_n+1}(c)}(P^{\pen{(n-qj_n)/q}+1})\\
&=\degl{x}\cdot\deg_{G_{n-qj_n+1}(c)}(P^{\pen{(n-qj_n)/q}+1})\\
&=\degl{x}\cdot d_{n-qj_n}(c),
\end{align*}
that is, $d_n$ belongs to 
\[\degl x\cdot\cD(c)\sub\degl x\cdot\bigcup_{c\in\cC(x)}\cD(c).\]
Now (\ref{item:unionlemmasigma}) follows.
\end{proof}

\begin{lemma}\label{lemma:rangeunboundedrec}
Let $P\in K[z]$ be a tame polynomial of degree $d\ge2$ and $x\in\JP.$ The range of dynamical degrees $\cD(x)$ is unbounded if and only if
there exists a recurrent critical element $c$ contained in the $\w$-limit of $x.$
\end{lemma}

To prove the previous lemma we need the following definition.

\begin{definition}\label{def:lcg}
For $x_1\ne x_2\in\JP$ the \emph{largest common geometric level between $x_1$ and $x_2$} is defined by
\[\lcg(x_1,x_2)\defi\max\{j\ge0\mid G_j(x_1)=G_j(x_2)\}.\] If $x_1=x_2$ we define $\lcg(x_1,x_2)$ as $+\infty.$
\end{definition}

\begin{proof}[Proof of Lemma~\ref{lemma:rangeunboundedrec}]
Suppose that $\cD(x)$ unbounded. By Lemma~\ref{lemma:sigma}~(\ref{item:unionlemmasigma}) it follows that there exists a critical element $c_1$
in $\barra{\cO^+(P(x))}$ such that $\cD(c_1)$ is also unbounded. Hence, we can find a critical element $c_2\in\barra{\cO^+(P(c_1))}$ with
$\cD(c_2)$ also unbounded. Recursively we obtain a sequence $(c_n)\sub\barra{\cO^+(P(x))}$ of critical elements with $\cD(c_n)$ unbounded and
$c_{n+1}\in\barra{\cO^+(P(c_n))},$ that is $\barra{\cO^+(P(c_{n+1}))}\sub\barra{\cO^+(P(c_n))}.$ By Proposition~\ref{prop:finitecritical} there are at most
$d-2$
critical elements of $P$ contained in $\JP.$ It follows that there exists $N\in\N$ such that $\barra{\cO^+(P(c_N))}$ coincides with
$\barra{\cO^+(P(c_{N+1}))}.$
That is, $c_{N+1}\in\barra{\cO^+(P(c_{N+1}))}.$ This means that $c_{N+1}$ is a recurrent critical element in $\barra{\cO^+(P(x))}.$ Hence, $c_{N+1}$
belongs to $\w(x).$

Conversely. Suppose that there exits a recurrent critical element $c$ in the $\w$-limit of $x.$

We split the proof in two parts. First we show that $\cD(c)$ is unbounded and then we prove that the range of dynamical degrees $\cD(x)$ is unbounded.

Suppose that $P$ has $q$ generators for its geometric sequences and denote by $(G_n)$ the geometric sequence of $c.$ The element $c$ is recurrent,
thus we can choose $n_1$ such that the largest common geometric level between $P^{n_1}(c)$ and $c$ is larger than $1,$ that is,
$\lcg(P^{n_1}(c),c)>1$ (see Definition~\ref{def:lcg}).
Hence, we have
\[\deg_{G_{1+qn_1}}(P^{\pen{(1+qn_1)/q}+1})=\deg_{G_{1+qn_1}}(P^{n_1})\cdot\deg_{G_1}(P^{\pen{(1+qn_1)/q}+1-n_1})\ge2^2.\]

Now we can pick $n_2>n_1$ such that $\lcg(P^{n_2}(c),c)>1+qn_1.$ Therefore
\begin{align*}
d_{1+qn_1+qn_2}(c)&=\deg_{G_{1+qn_1+qn_2}}(P^{\pen{(1+qn_1+qn_2)/q}+1})\\
&=\deg_{G_{1+qn_1+qn_2}}(P^{n_2})\cdot\deg_{G_{1+qn_1}}(P^{\pen{(1+qn_1+qn_2)/q}+1-n_2})\\
&=\deg_{G_{1+qn_1+qn_2}}(P^{n_2})\cdot\deg_{G_{1+qn_1}}(P^{\pen{(1+qn_1)/q}+1})\\
&\ge2^2\cdot\deg_{G_{1+qn_1+qn_2}}(P^{n_2})\\
&\ge2^3.
\end{align*}

Recursively, we can find a increasing sequence $(n_j)$ of natural numbers such that, for $\sigma_j=1+qn_1+\cdots+qn_j$
\[\deg_{G_{\sigma_j}}(P^{i_{\sigma_j}+1})\ge2^{j+1}.\]

It follows that $\cD(c)$ is an unbounded set.

For each $n\ge0$ there exists at least one $m_n\in\N$ such that $G_{m_n}(x)$ is mapped on $G_n$ by some iterate $P^{j_n}$ of $P.$ Then
\begin{align*}
d_{m_n}(x)&=\deg_{G_{m_n}(x)}(P^{\pen{m_n/q}+1})\\
&=\deg_{G_{m_n}(x)}(P^{j_n})\cdot\deg_{G_n(c)}(P^{\pen{m_n/q}+1-j_n})\\
&=\deg_{G_{m_n}(x)}(P^{j_n})\cdot\deg_{G_n(c)}(P^{\pen{n/q}+1})\\
&\ge\degl x\cdot\deg_{G_n(x)}(P^{\pen{n/q}+1})\\
&=\degl x\cdot d_n(c)
\end{align*}
Hence, the range of dynamical degrees $\cD(x)$ is unbounded.
\end{proof}

\begin{lemma}\label{lemma:snxunbimpdxunb}
Let $P\in k[z]$ be a nonsimple polynomial of degree $\ge2.$ Suppose that there exists $x\in\JP$ such that
$\{s_n(x)\mid t_n\ne t_{n+1}\}$ is unbounded. Then the range of dynamical degrees around $x$ is unbounded.
\end{lemma}
\begin{proof}
If $\{s_n(x)\mid t_n\ne t_{n+1}\}$ is unbounded we have that
\[\log(\diam(G_j))=\log(\diam(L_0))-\dH(G_j,L_0)\]
belongs to $k_j^{-1}\N$ for $(k_j)$ an unbounded sequence, then, in view of Lemma~\ref{lemma:sumafea},
\[\liminf_{n\to\infty}[d_n(x)]^{-1}=0.\]
Hence, the range of dynamical degrees $\cD(x)$ is unbounded.
\end{proof}

\begin{lemma}\label{lamma:recimptrans}
Let $P\in k[z]$ be a tame polynomial of degree $\ge2.$ Suppose that there exists a recurrent critical element
$c\in\JP.$ Then, there exists a sequence $(y_n)\sub\JP$ such that $\displaystyle{\lim_{n\to+\infty} y_n=c}$ and $\displaystyle{\lim_{n\to+\infty}\delta(y_n)=+\infty}.$
\end{lemma}

To prove the previous lemma it is convenient to use $p$-adic absolute value notation. Recall that we are supposing that
\[\log(\abs{k^\times})=\Z\quad\textnormal{and}\quad\log(\abs{K^\times})=\Q.\]
Given a prime number $p,$ the absolute
value $\absp\cdot$ is defined as follows. For $a\in\Z$ we put $\absp a=p^{-n},$ where, $a=p^{n}m$ and $p$ does not divides $m.$ For a rational 
number $\displaystyle{\frac{a}{b}}$ we put $\displaystyle{\absp{\frac{a}{b}}=\frac{\absp a}{\absp b}}.$ Recall that $\absp\cdot$ is a non-Archimedean
absolute value.

For convenience we record the following general fact about non-Archimedean absolute values as a remark.

\begin{remark}\label{rem:mayorgana}
Let $r_1,r_2\in\Q$ such that $\absp{r_1}<\absp{r_2},$ then $\absp{r_1+r_2}=\absp{r_2}.$
\end{remark} 

\begin{proof}[Proof of Lemma~\ref{lamma:recimptrans}]
In view of Lemma~\ref{lemma:rangeunboundedrec} we have that the range of dynamical degrees $\cD(c)$ is an unbounded set.
For each $n\ge0$ the dynamical degree $d_n=d_n(c)$ around $c$ is a product of $\pen{n/q}+1$ numbers smaller than $d.$ Hence, there exists a prime
number $p\le d$ such that arbitrarily large powers of $p$ divide elements in the range of dynamical degrees around $c,$ that is, $(\absp{d_n^{-1}})$ 
is an unbounded sequence.

If $I_{j_n}$ denotes the segment such that $P^{\pen{n/q}+1}(]G_{n},G_{n-1}[)=I_{j_n}$ and
\[\alpha_n=\absp{d_n^{-1}\cdot\abs{I_{j_n}}},\] 
it follows that the sequence $(\alpha_n)$ is also unbounded.

Let $m_1\ge1$ be the smallest integer such that $\max\{1,\absp{\log(\diam(L_0))}\}<\absp{\alpha_{m_1}}.$ Let $1<e_1$ such that
$\absp{\alpha_{m_1}}=p^{e_1}.$ 

Recursively, we define $m_n$ as the smallest integer larger than $m_{n-1}$ such that $p^{e_{n-1}}<\absp{\alpha_{m_n}}.$ We define $e_n$ as the
integer 
satisfying $\absp{\alpha_{m_n}}=p^{e_n}.$ Note that both $(m_n)$ and $(e_n)$ are increasing sequences of integers.

To construct the sequence $(y_n)$ we need to prove first that $\{s_{m_n}(c)\mid n\in\N\}$ is an unbounded set. We have that
\begin{align*}
-\log(\diam(G_{m_n}))&=-\log(\diam(L_0))+\sum_{\ell=1}^{m_n}d_\ell^{-1}\cdot\abs{I_{j_\ell}}\\
&=-\log(\diam(L_0))+\sum_{\ell=1}^{m_n}\alpha_\ell.
\end{align*}
Using Remark~\ref{rem:mayorgana} and the choice of $m_n$ we have that
\begin{align*}
\absp{\log(\diam(G_{m_n}))}&=\absp{-\log(\diam(L_0))+\sum_{\ell=1}^{m_n}\alpha_\ell}\\
&=\absp{\alpha_{m_n}}=p^{e_n}.
\end{align*}

Hence, the largest power of $p$ dividing the denominator, say $b_n,$ of $\log(\diam(G_{m_n}))$ is $p^{e_n}.$ Then
\[s_{m_n}=\frac{\lcm(m,b_j)}{m}\ge\frac{\lcm(m,p^{e_n})}{m},\]
by definition of $s_n$ (see Definition~\ref{def:sn}). In particular, $\lim s_{m_n}=+\infty,$ since $(e_n)$ is an increasing
sequence of integers.

For all $n\ge q\cN(c)+q$ there exists a noncritical element $y_n\in\JP$ such that $G_{m_n}(y_n)=G_{m_n}$ and $G_{m_n+1}(y_n)\ne G_{m_n+1},$ see  
Note that, $\lcg(y_n,c)=m_n\to+\infty$ as $n\to+\infty.$ Then the sequence $(y_n)$
converges to $c.$ In view of this, we only need prove that $\lim\delta(y_n)=+\infty$

The injectivity time of $G_{m_n+1}(y_n)$ is larger than the injectivity time of $G_{m_n}(y_n),$ since $G_{m_n}(y_n)$ is a critical element and
$G_{m_n+1}$
is not. Hence, applying Proposition~\ref{prop:condicion} and Lemma~\ref{lemma:justi}~(\ref{item:leclemmajusti}) we
obtain a bound, from below, for the algebraic degree of $y_n.$ More precisely,
\[\delta(y_n)\ge\delta(G_{m_n+1}(y_n))\ge\max\{1,s_{m_n}(c)\}=s_{m_n}.\]
That is, $\lim\delta(y_n)=+\infty.$
\end{proof}

Now we can give the proof of Theorem~\ref{theo:juliaalg}.

\begin{proof}[Proof of the Theorem~\ref{theo:juliaalg}]
We prove that in the presence of a recurrent critical element, $\JP$ is not an algebraic set over $k.$ Suppose that there
exists a recurrent critical element $c$ contained in $\JP.$ Following Lemma~\ref{lamma:recimptrans}, there exists a sequence
$(y_n)\sub\JP$ such that $\lim\delta(y_n)=+\infty.$ Applying Corollary~\ref{cor:unotodo} we obtain $x$ in the Julia set of $P$
with $\delta(x)=+\infty.$

In particular, if $\JP$ is algebraic over $k,$ then there are not critical periodic elements in $\JP.$ Following
Proposition~\ref{prop:algnowand} we obtain that $\JP$ is contained in $K.$ In Proposition~\ref{prop:Delta} we showed the
existence of (the smallest) $\Delta\in\N$ with $\delta(x)\le\Delta$ for all $x\in\JP.$
\end{proof}

\subsection{The completion of the field of formal Puiseux series}\label{subsec:lau}

Let $F$ be an algebraically  closed field of characteristic 0. We denote by $\flau$ the \emph{field of formal Laurent series} with
coefficients in $F.$ For a nonzero element
\[z=\sum_{j\ge j_0} a_j\tau^j\in\flau\]
we define $\ord(z)=\min\{a_j\mid j\ne0\}$ and $\abs z=e^{-\ord(z)}.$
Observe that $\abs\cdot$ is a non-Archimedean absolute value and that $\flau$ is complete with respect to $\abs\cdot$ but not algebraically closed.
An algebraic closure of $\flau$ is the \emph{field of formal Puiseux series $\fpui$} with coefficients in $F$ (e.g. 
chapter~IV, Theorem~3.1 in \cite{algebraic_walker}). More precisely, $\fpui$ is the direct limit of the fields $\flau(\tau^{1/m})=\flaum$ for $m\in\N,$
with the obvious inclusions, that is, $F(\!(\tau^{1/{m_1}})\!)\sub F(\!(\tau^{1/{m_2}})\!)$ if and only if $m_1$ divides to $m_2.$

Therefore an element in $\fpui$ has the form
\[z=\sum_{j\ge j_0} a_j\tau^{j/m}\] for some $m\in \N.$ The unique extension of $\abs\cdot$ to $\fpui$ (also denoted by $\abs\cdot$) is completely
determined by $\abs{\tau^{1/m}}=e^{-1/m}.$

Note that the degree of the field $\flaum$ over $\flau$ is precisely $m.$ Moreover, $\flaum$ is the unique field extension of $\flau$ of degree $m.$

We denote by $\L_F$ the completion of $\fpui$ with respect to $\abs\cdot.$ Every $z\in\L_F$ can be represented as
\[\sum_{j\ge 0} a_{j}\tau^{\lambda_j},\]
where $(\lambda_j)$ is an increasing sequence of rational numbers tending to $+\infty.$ In this case $\abs z=e^{-\ord(z)}$ where
$\ord(z)=\min\{\lambda_j\mid a_j\ne0\}.$ The field $\L_F$ is also algebraically closed since it is  the completion of an algebraically
closed non-Archimedean field.

The following Proposition is a complement of Proposition~\ref{prop:condicion}. The result is not true in general, it depends on the structure of the finite extensions of $\flau.$

\begin{proposition}\label{prop:lcmcondicion}
Let $P\in\fpui[z]$ be a nonsimple polynomial of degree $\ge2.$ Consider $x$ in the Julia set of $P$ and denote by
$(G_n)$ the geometric sequence of $x.$ For each $n\ge-1$ let $s_n=s_n(x).$ Then
\begin{equation}\label{eq:lcmpropocondicion}
\delta(G_{n+1})\le\lcm\{s_j\mid t_j\ne t_{j+1},-1\le j\le n\},
\end{equation}
where $\lcm$ denotes the least common multiple.
\end{proposition}
\begin{proof}
We proceed by induction. If $n=-1$ we have that
\[\delta(G_0)=1\le\lcm\{s_j\mid t_j\ne t_{j+1},-1\le j\le-1\}.\]

Suppose that (\ref{eq:lcmpropocondicion}) holds for $n-1.$ If $t_n=t_{n+1}$ then
\[\{s_j\mid t_j\ne t_{j+1},-1\le j\le n-1\}=\{s_j\mid t_j\ne t_{j+1},-1\le j\le n\}.\]
By Corollary~\ref{cor:tndifferent} we have that $\delta(G_n)=\delta(G_{n+1}),$ therefore (\ref{eq:lcmpropocondicion}) holds for $n.$

Suppose that $t_{n}\ne t_{n+1}$ and let $\cD_n$ be the direction in $T_{G_n}$ containing
$G_{n+1}.$ We have that, there exists $0\ne a\in F$ such that $u_n+a\tau^{b_n/a_n}$ belongs to $\cD_n\cap\L_F$ where
$b_n/a_n=\abs{\log\diam(G_n)}_\R$ with $(a_n,b_n)=1.$

In view of the structure of the algebraic extensions of the field of formal Laurent series (see subsection~\ref{subsec:lau}) we have that
$\delta(u_n+a\tau^{b_n/a_n})\le\lcm\{\delta(\tau^{b_n/a_n}),\delta(u_n)\}$ and $\delta(\tau^{b_n/a_n})=\lcm\{m,a_n\}\cdot m^{-1}=s_n.$
It follows, by Lemma~\ref{lemma:degree}, that
\begin{align*}
\delta(u_{n+1})&\le\delta(u_n+a\tau^{b_n/a_n})\\
&\le\lcm\{\delta(\tau^{b_n/a_n}),\delta(u_n)\}\\
&=\lcm\{\delta(u_n),s_n\}
\end{align*}
Hence, we have that $\delta(G_{n+1})\le\lcm\{s_n,\delta(G_n)\}.$

Applying the inductive hypothesis it follows that
\begin{align*}
\delta(G_{n+1})&\le\lcm\{s_n,\delta(G_n)\}\\
&\le\lcm\big\{s_n,\lcm\{s_j\mid t_j\ne t_{j+1},-1\le j\le n-1\}\big\}\\
&\le\lcm\{s_j\mid t_j\ne t_{j+1},-1\le j\le n\}.
\end{align*}
This proves the proposition.
\end{proof}

\begin{proof}[Proof of Theorem~\ref{theo:juliaalgiff}]
From Theorem~\ref{theo:juliaalg} we have that there are not recurrent critical element in $\JP$ and that there exists a smallest
$\Delta\in\N$ such that $\delta(x)\le\Delta$ for all $x\in\JP.$ In particular, in view of the structure of the subfields of $\fpui,$
it follows that $\JP$ is contained in the unique extension of $\flau$ with degree $\Delta.$

Conversely. Suppose that there exists $x\in\JP$ with $\delta(x)=+\infty.$
By Proposition~\ref{prop:lcmcondicion} it follows that $\lcm\{s_n(x)\mid t_n\ne t_{n+1}\}$ is an unbounded set,
hence $\{s_n(x)\mid t_n\ne t_{n+1}\}$ is also an unbounded set and therefore the range of
dynamical degrees $\cD(x)$ is also unbounded (see Lemma~\ref{lemma:snxunbimpdxunb}). By
Lemma~\ref{lemma:rangeunboundedrec}, there exists a critical element $c\in\JP$ with its range of dynamical
degrees $\cD(c)$ unbounded. Hence, there exist a recurrent critical element in $\JP.$
\end{proof}

\section{Polynomial with Coefficients in $K$}\label{sec:polele}

In the previous section we proved Corollary~\ref{cor:nowandering} for polynomials with algebraic coefficients. To
prove Theorem~\ref{theo:juliaberk} for polynomials with coefficients in $K$ we use a perturbation argument.

The key to perturb the coefficients of $P$ preserving a suitable orbit is the following proposition.

\begin{proposition}\label{prop:diam}
Let $P\in K[z]$ be a tame polynomial of degree $d\ge 2$ and consider $x\in\JP.$ Then
\[\dH(y,L_0)\le d^{d-1}\cdot\dH(x,L_0),\]
for all $y$ in $\w(x).$
\end{proposition}

\subsection{Key Lemma and Proof of Proposition~\ref{prop:diam}}

To prove Proposition~\ref{prop:diam} we need to compare the distance from $\dH(x,L_0)$ with $\dH(y,L_0).$ To do this, we need to introduce
the concepts of \emph{level} and \emph{time sequences.} Recall that $\lcg(x,y)$ denotes the largest common geometric level between $x$ and $y,$ see Definition~\ref{def:lcg}.  

Let $P\in K[z]$ be a nonsimple polynomial of degree $\ge2$ with $q$ generators for its geometric sequences.

We define $k_0=0$ and $\ell_0=\lcg(x,y).$ The point $y$ belongs to $\w(x),$ hence
\[\limsup_{j\to+\infty}\lcg(P^j(x),y)=+\infty.\] Let
\[k_1=\min\{j\in\N\mid\lcg(P^j(x),y)>\ell_0\}\]
and let $\ell_1=\lcg(P^{k_1}(x),y).$

\begin{figure}
\begin{center}
\includegraphics{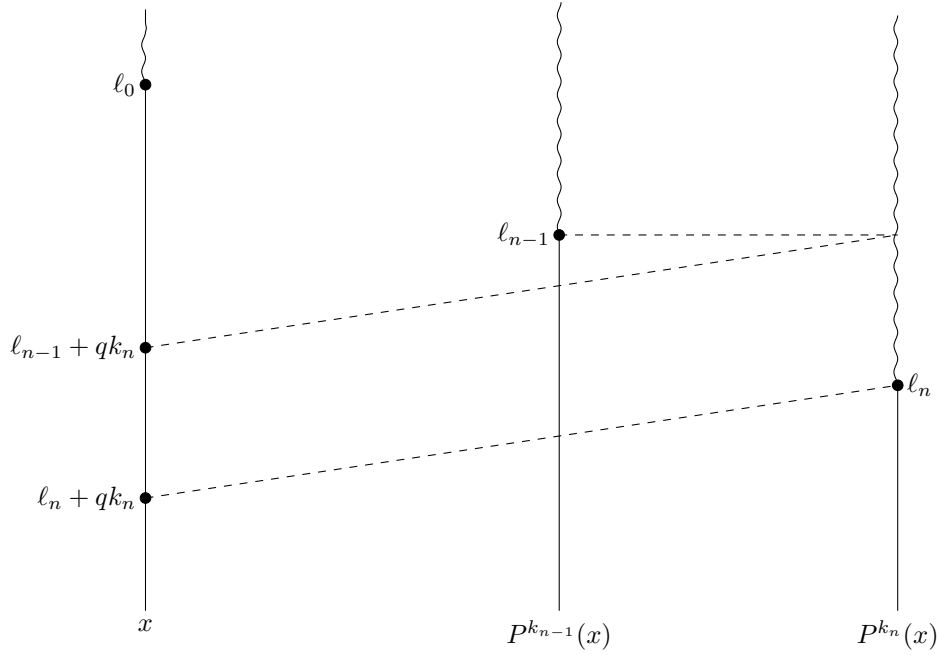}
\caption{Level and time sequences}
\label{fig:leveltime}
\end{center}
\end{figure}

Recursively, we define $k_n=\min\{j\in\N\mid\lcg(P^j(x),y)>\ell_{n-1}\}$ and
\[\ell_n=\lcg(P^{k_n}(x),y)\] for all $n\ge2.$ See Figure~\ref{fig:leveltime}.

\begin{definition}\label{def:lnkn}
We say that $(\ell_n)_{n\ge0}$ is the \emph{level sequence from $x$ to $y$} and that $(k_n)_{n\ge0}$ is the \emph{time sequence from $x$ to $y.$}
\end{definition}

Note that
\[P^{k_n}([G_{\ell_n+qk_n}(x),G_{\ell_{n-1}+qk_n}(x)])=[G_{\ell_n}(y),G_{\ell_{n-1}}(y)].\]

\begin{remark}\label{rem:knmin}
Observe that $k_n$ is the smallest integer such that
\[P^{k_n}([x,L_0])\cap[y,G_{\ell_{n-1}}(y)[\,\ne\vacio.\]
\end{remark}

\begin{lemma}\label{lemma:gradoacotado}
Let $P\in K[z]$ be a tame polynomial of degree $d\ge 2$ with $q$ generators for its geometric sequences. Let $x\in\JP$ and consider $y$ in
the $\w$-limit of $x,$ we denote by $(\ell_n)$ (resp. $(k_n)$) the level (resp. time) sequence from $x$ to $y.$ Then, for any $n\ge1,$
\[\deg_{G_j(x)}(P^{k_n})\le d^{d-1}\]
for all the elements $G_j(x)$ contained in the segment $[x,G_{\ell_{n-1}+qk_n}(x)[.$
\end{lemma}
\begin{proof}
Let $n>1$ and let $G_j(x)\in[x,G_{\ell_{n-1}+qk_n}(x)[$ be a geometric point. Note that the elements in the set
\[O_n(G_j(x))=\{G_j(x),P(G_j(x)),\dots,P^{k_n-1}(G_j(x))\}\]
are pairwise incomparable. In fact, if we suppose that there exist $0\le i$ and $1\le\ell$ such that $1\le i+\ell<k_n$ and
$P^{i}(G_j(x))\lec P^{i+\ell}(G_j(x)),$ 
it follows that
\[G_{j-qi}(P^i(x))\lec G_{j-qi-q\ell}(P^{i+\ell}(x)).\]
Hence 
\[G_{j-qi-q\ell}(P^i(x))=G_{j-qi-q\ell}(P^{i+\ell}(x)).\]
Then we have that
\begin{align*}
P^{k_n-\ell}(G_{j-q\ell}(x))&=P^{k_n-i-\ell}(P^i(G_{j-q\ell}(x)))\\
&=P^{k_n-i-\ell}(G_{j-qi-q\ell}(P^i(x))\\
&=P^{k_n-i-\ell}(G_{j-qi-q\ell}(P^{i+\ell}(x))\\
&=G_{j-qi-q\ell-q(k_n-i-\ell)}(P^{i+\ell+k_n-i-\ell}(x))\\
&=G_{j-qk_n}(P^{k_n}(x))\\
&=P^{k_n}(G_j(x)).
\end{align*}

Since $P^{k_n}(G_j(x))\in[y,G_{\ell_{n-1}}(y)[,$ we have $P^{k_n-\ell}([x,L_0])\cap[y,G_{\ell_{n-1}}(y)[\,\ne\vacio,$ 
which contradicts Remark~\ref{rem:knmin}.

Hence, we obtain that the set $O_n(G_j(x))$ contains at most $d-1$ critical elements. In consequence, the local degree of $P^{k_n}$ at $G_j(x)$
is bounded, from above, by $d^{d-1}.$
\end{proof}

\begin{proof}[Proof of Proposition~\ref{prop:diam}]
Let $n\ge1$ and suppose that $P$ has $q$ generators for its geometric sequences. Let $x$ be in the Julia set of $P$ and consider $y\in\w(x).$ 
Denote by $(\ell_n)$ (resp. $(k_n)$) the level (resp. time) sequence from $x$ to $y.$ From the previous Lemma and Proposition~\ref{prop:enecalig}, the
local degree of $P^{k_n}$ is constant and smaller than $d^{d-1}$ in each segment of the form $]G_{j+1}(x),G_j(x)[$ contained
in $[G_{\ell_n+qk_n}(x),G_{\ell_{n-1}+qk_n}(x)[.$ By Lemma~\ref{lemma:distgrado} it follows that
\[\dH(P^{k_n}(G_{j+1}(x)),P^{k_n}(G_j(x)))\le d^{d-1}\cdot\dH(G_{j+1}(x),G_j(x)).\]

Applying this to all the segments $]G_{j+1}(x),G_j(x)[$ contained in
\[[G_{\ell_n+qk_n}(x),G_{\ell_{n-1}+qk_n}(x)[\] we obtain that
\[\dH(G_{\ell_n}(y),G_{\ell_{n-1}}(y))\le d^{d-1}\cdot\dH(G_{\ell_n+qk_n}(x),G_{\ell_{n-1}+qk_n}(x)).\]

Therefore, if we put $a=\dH(G_{\ell_0}(x),L_0)=\dH(G_{\ell_0}(y),L_0),$ we have
\begin{align*}
\dH(y,L_0)&=a+\sum_{n=1}^{+\infty}\dH(G_{\ell_{n}}(y),G_{\ell_{n-1}}(y))\\
&\le a+d^{d-1}\cdot\sum_{n=1}^{+\infty}\dH(G_{\ell_n+qk_n}(x),G_{\ell_{n-1}+qk_n}(x))\\
&\le a+d^{d-1}\cdot\sum_{n=\ell_0+qk_1}^{+\infty}\dH(G_{n+1}(x),G_{n}(x))\\
&\le a+\dH(G_{\ell_0+qk_1}(x),G_{\ell_0}(x))+d^{d-1}\cdot\dH(x,G_{\ell_0+qk_1}(x))\\
&\le d^{d-1}\cdot\dH(x,L_0)
\end{align*}
\end{proof}

\subsection{Recurrent Orbits in the Hyperbolic Space}

Let $P\in K[z]$ be a tame polynomial of degree $\ge2$ and consider $x\in\JP\setminus K.$ From Proposition~\ref{prop:hlrec} we know that
$\w(x)$ contains a critical element. Nevertheless, using Proposition~\ref{prop:diam} we can be more precise.

\begin{corollary}\label{cor:hpcritrec}
Let $P\in K[z]$ be a tame polynomial of degree $\ge2$ and consider $x\in\JP\setminus K.$ Then, $\w(x)\cap\HK$ contains a recurrent critical
element.
\end{corollary}
\begin{proof}
Let $x\in\JP\cap\HK.$ From Proposition~\ref{prop:hlrec} there exists a critical point $c_1$ in $\w(x).$ By Proposition~\ref{prop:diam} it follows that
\[\dH(c_1,L_0)\le d^{d-1}\cdot\dH(x,L_0)<+\infty.\]
Hence, $c_1$ belongs to $\HK.$ Recursively, we can find a sequence of critical elements $(c_n)\sub\JP\cap\HK$ such that $c_{n+1}\in\w(c_{n}).$ By
Proposition~\ref{prop:finitecritical} there exists $N$ such that $c_N\in\w(c_N).$ That is, $c_N$ is a recurrent critical element in $\w(x)\cap\HK.$ 
\end{proof}

By corollary above, we need to study the recurrent critical elements in $\JP\cap\HK.$ 

In the case that $\w(x)$ contains a periodic critical orbit we will show, Proposition~\ref{prop:pointsperiodic}, that $x$ is a preperiodic point.
In Proposition~\ref{prop:critrat}, we will prove that the unique recurrent critical orbits in $\JP\cap\HK$ are the periodic critical orbit.

\begin{proposition}\label{prop:pointsperiodic}
Let $P\in K[z]$ be a tame polynomial of degree $\ge2$ and consider $x\in\JP.$ If the $\w$-limit of $x$ contains a periodic critical element,
then $x$ is preperiodic or $x$ is a classical point.
\end{proposition}

\begin{proposition}\label{prop:critrat}
Let $P\in K[z]$ be a tame polynomial of degree $\ge2.$ Then the critical elements of $P$ contained in $\JP\cap\HK$ are preperiodic critical elements of
type II.
\end{proposition}

To show Proposition~\ref{prop:critrat} we use a perturbation argument. To prove Proposition~\ref{prop:pointsperiodic} we need the following lemma
about geometric sequences.

\begin{lemma}\label{lemma:previous}
Let $P\in K[z]$ be a tame polynomial of degree $\ge2$ with $q$ generators for its geometric sequences. Consider $c\in\Crit(P)\cap\JP$ a critical
element  which is fixed by $P.$ Let $x$ be a Julia point such that $G_n(x)=G_n(c)$ and $G_{n+1}(x)$ is noncritical for some $n>q+q\cN(c).$
Then $P(G_{n+1}(x))$ is not critical.
\end{lemma}
\begin{proof}
In view of the properties of $\cN(c)$ (see Proposition~\ref{cor:M}) we know that if $P(G_{n+1}(x))$ is critical, then it coincides with
$G_{n+1-q}(c)=G_{n+1-q}(P(c)).$
Therefore \[\degl{G_n(c)}\ge1+\degl{G_{n+1}(c)},\]
from Remark~\ref{rem:RH}, which contradicts the definition of $\cN(c),$ since the local degrees of $P$ at $G_n(c)$ and at $G_{n+1}(c)$ coincide.
\end{proof}

\begin{proof}[Proof of Proposition~\ref{prop:pointsperiodic}]
We assume that $P$ is a tame polynomial with $q$ generators for its geometric sequences.

Suppose that there exists a nonpreperiodic element $x\in\JP$ such that the $\w$-limit of $x$ contains
a periodic critical element $c.$ Passing to an iterate we can suppose that $c$ is a fixed point.

Let $(\ell_n)$ be the level sequence from $x$ to $c$ and $(k_n)$ be the time sequence from $x$ to $c.$

Since $1+\ell_n+qk_n>\ell_{n-1}+qk_n$ we can use Lemma~\ref{lemma:gradoacotado} to conclude that
\begin{equation}\label{eq:dd-1}
\deg_{G_{1+\ell_n+qk_n}(x)}(P^{k_n})\le d^{d-1}
\end{equation}
for all $n\ge1.$

Consider $\cN=\max\{\cN(x),\cN(c)\}.$ If $\ell_n>q+q\cN$ the points
\[G_{1+\ell_n}(P^{k_n}(x))=P^{k_n}(G_{1+\ell_n+qk_n}(x))\textnormal{ and }G_{\ell_n}(P^{k_n}(x))=P^{k_n}(L_{\ell_n+qk_n}(x))\]
satisfy the hypothesis of Lemma~\ref{lemma:previous}, since $G_{\ell_n}(P^{k_n}(x))=G_{\ell_n}(c)$ and 
\[G_{1+\ell_n}(P^{k_n}(x))=P^{k_n}(G_{1+\ell_n+qk_n}(x))\]
is not a critical element, by definition of $(\ell_n)$ and $(k_n)$ (see Definition~\ref{def:lnkn}). Hence
\[P(G_{1+\ell_n}(P^{k_n}(x)))\]
is not a critical element, that is $\deg_{G_{1+\ell}(P^{k_n}(x))}(P^{2})=1.$ If $\ell_n-q>q-q\cN$ the points
\[G_{1+\ell_n-q}(P^{k_n+1}(x))=P(G_{1+\ell_n}(P^{k_n}(x)))\]
and
\[G_{\ell_n-q}(P^{k_n+1}(x))=P(G_{\ell_n}(P^{k_n}(x)))=P(G_{\ell_n}(c))=G_{\ell_n-q}(c)\]
satisfy the hypothesis of Lemma~\ref{lemma:previous} by the definition of the level and time sequences. Hence
\[P(G_{1+\ell_n-q}(P^{k_n+1}(x)))=P^2(G_{1+\ell_n}(P^{k_n}(x)))\]
is not a critical element, that is $\deg_{G_{1+\ell}(P^{k_n}(x))}(P^{3})=1$. Applying Lemma~\ref{lemma:previous} recursively we obtain that the elements in the
\[\{P(G_{1+\ell_n}(P^{k_n}(x))),P^2(G_{1+\ell_n}(P^{k_n}(x))),\dots,P^{e_n-1}(G_{1+\ell_n}(P^{k_n}(x)))\}\]
are not critical, where $e_n$ is the smallest positive integer such that 
\[\ell_n-q(e_n-2)\le q+q\cN.\]
In particular, it follows that
\begin{equation}\label{eq:d=1}
\deg_{G_{1+\ell}(P^{k_n}(x))}(P^{e_n})=\deg_{G_{1+\ell_n+qk_n}(x)}(P^{k_n+e_n})=1.
\end{equation}

If we have that $\ell_n-q(e_n-2)\le q+q\cN$ it follows that
\[\ell_n-qe_n\le -q+q\cN\le q\cN,\]
since $1\le q.$ Then
\begin{align*}
\deg_{G_{1+\ell_n-qe_n}}(P^{\pen{\ell_n/q}-e_n+1})&=\deg_{G_{1+\ell_n-qe_n}}(P^{\pen{(\ell_n-qe_n)/q}+1})\\
&\le\deg_{G_{1+\ell_n-qe_n}}(P^{\pen{q\cN/q}+1})\\
&=\deg_{G_{1+\ell_n-qe_n}}(P^{\cN+1})\\
&\le d^{\cN+1}
\end{align*}
that is
\begin{equation}\label{eq:en}
\deg_{G_{1+\ell_n-qe_n}}(P^{\pen{\ell_n+qe_n/q}+1})=\deg_{G_{\ell_n-qe_n}}(P^{\pen{\ell_n/q}+e_n+1})\le d^{\cN+1}
\end{equation}
for any geometric point of level $1+\ell_n-qe_n.$

If
\[\Delta_n=\deg_{G_{1+\ell_n+qk_n}(x)}(P^{\pen{(\ell_n+qk_n)/q}+1})\]
we have that
\begin{align*}
\Delta_n&=\deg_{G_{1+\ell_n+qk_n}(x)}(P^{\pen{\ell_n/q}+k_n+1})\\
&=\deg_{G_{1+\ell_n+qk_n}(x)}(P^{k_n})\cdot\deg_{P^{k_n}(G_{1+\ell_n+qk_n}(x))}(P^{e_n})\cdot\deg_{P^{k_n+e_n}(G_{1+\ell_n+qk_n}(x))}(P^{\pen{\ell_n/q}+e_n+1})\\
&=\deg_{G_{1+\ell_n+qk_n}(x)}(P^{k_n})\cdot\deg_{G_{1+\ell_n}(P^{k_n}(x))}(P^{e_n})\cdot\deg_{G_{1+\ell_n-qe_n}(P^{k_n+e_n}(x))}(P^{\pen{\ell_n/q}+e_n+1})\\
\end{align*}
by (\ref{eq:dd-1}), (\ref{eq:d=1}) and (\ref{eq:en}) we conclude that
\begin{equation}\label{eq:alfin}
\deg_{G_{1+\ell_n+qk_n}(x)}(P^{\pen{(1+\ell_n+qk_n)/q}+1})\le d^{d-1}\cdot 1\cdot d^{\cN+1}=d^{d+\cN}
\end{equation}
for all $n\ge M$ where $M\in\N$ is such that $\ell_M>q+q\cN.$

Following Proposition~\ref{prop:enecalig} and (\ref{eq:alfin}) above, we have that $\deg_y(P^{\pen{(\ell_n+qk_n)/q}+1})$ is constant and bounded by $d^{d+\cN}$ for any point $y$ in the segment
\[[G_{1+\ell_n+qk_n}(x),G_{\ell_n+qk_n}(x)[.\]

Since $P^{\pen{(\ell_n+qk_n)/q}+1}([G_{1+\ell_n+qk_n}(x),G_{\ell_n+qk_n}(x)])$ belong to $\{I_0,I_1,\dots,I_{q-1}\}$ (see discussion below Definition~\ref{def:geo}) there exists
a segment $I\in\{I_0,I_1,\dots,I_{q-1}\}$ and an strictly increasing sequence $(n_j)$ of integers larger that $M$ such that
\[P^{\pen{(\ell_{n_j}+qk_{n_j})/q}+1}([G_{1+\ell_{n_j}+qk_{n_j}}(x),G_{\ell_{n_j}+qk_{n_j}}(x)])=I,\]
for all $j\ge1.$

Therefore, following Lemma~\ref{lemma:distgrado} and using (\ref{eq:alfin}), we have that
\begin{align*}
\dH(x,L_0)&\ge\sum_{j=1}^{+\infty}\dH(G_{1+\ell_{n_j}+qk_{n_j}}(x),G_{\ell_{n_j}+qk_{n_j}}(x))\\
&\ge\sum_{j=1}^{+\infty}\frac{\abs{I}}{d^{d+\cN}}\\
&=+\infty,
\end{align*}
where $\abs{I}$ denotes the length of the segment $I.$ It follows, $x$ is a classical point.
\end{proof}

Recall that we are assuming that $k\sub K$ is a complete field such that the restriction of $\abs\cdot$ to $k$ is a discrete absolute value and that
\[\{z\in K\mid [k(z),k]<+\infty\}\]
is a dense subset of $K.$

In order to prove Proposition~\ref{prop:critrat} we need a proposition to relate the dynamics of a polynomial
\[P(z)=a_dz^d+\cdots+a_1z+a_0\]
in $ K[z]$ with the dynamics of a small perturbation $Q$ of $P,$ of degree $d$ and with coefficients in a finite extension of $k.$

\begin{proposition}\label{prop:aprox}
Let $P\in K[z]$ be a tame polynomial of degree $d\ge2$ and let $0<R<\diam(L_0).$ Then, a tame polynomial $Q\in\calg k[z]$ such that $P(x)=Q(x)$ and $\degl x=\deg_x(Q)$ for all $x$ in
\[\barra{D_0}\cap\{\diam(y)>R\}\cap P^{-1}(\{\diam(y)>R\}).\] 
\end{proposition}

To prove the previous Proposition we need the following straightforward fact that we establish as a lemma

\begin{lemma}\label{lemma:cuenta}
Let $\eta>0.$ For each $n\in\N$ consider $w_1,\dots,w_n,w_1',\dots,w_n'\in K$ such that $\abs{w_j-w_j'}<\eta$ and $\abs{w_j}=\abs{w_j'}$ for all $1\le j\le n.$ Then
\[\abs{w_1w_2\cdots w_n-w_1'w_2'\cdots w_n'}\le(\max\{\abs{w_1},\abs{w_2},\dots,\abs{w_n}\})^{n-1}\cdot\eta.\]
\end{lemma}

\begin{proof}[Proof of Proposition~\ref{prop:aprox}]
Let $P(z)=a_dz^d+a_{d-1}z^{d-1}+\cdots+a_1z+a_0\in K[z]$ be a tame polynomial. Since $P$ is a tame polynomial we have that
\[\Crit(P)=\bigcup_{j=1}^{d-1}[w_j,\infty[,\]
where $P'(z)=a_d\!\cdot\! d\!\cdot\!(z-w_1)(z-w_2)\cdots(z-w_{d-1}).$ Note that
\[a_j=\frac{a_d\!\cdot\!d}{j}\cdot\sum_{1\le i_1<i_2<\cdots<i_{d-j}\le{d-1}}w_{i_1}w_{i_2}\cdots w_{i_{d-j}}\]
for all $1\le j\le d-1.$

Let
\begin{equation}\label{eq:C}
C=\max_{1\le j\le d-1}\Big\{\abs{j^{-1}}\Big\}\cdot\max\bigg\{1,\abs{a_d}^d,\Big(\max_{1\le j\le d-1}\{\abs{w_j}\}\Big)^d\bigg\}.
\end{equation}

Put $D_0=B_{r_P}(z_0)$ and let
\begin{equation}\label{eq:r0}
R_0=\max\{r_P,\abs{z_0}\}
\end{equation}

Consider
\begin{equation}\label{eq:epsilon}
0<\e<\min\bigg\{\min_{w_i\ne w_j}\{\abs{w_i-w_j}\},R,\frac{R}{C\cdot\max\{1,R_0^d\}}\bigg\}.
\end{equation}

Since $\calg k$ is dense in $K,$ for all $1\le j\le d-1$ we can pick $w_j'\in\calg k$ such that $\abs{w_j-w_j'}<\e,\,\abs{w_j}=\abs{w_j'}$ and $w_j'=w_i'$ if $w_j=w_i.$ Moreover, pick $b_d,b_0\in\calg k$ such that
$\abs{b_j-a_j}<\e$ and $\abs{a_j}=\abs{b_j}$ for $j=0,d.$ Let
\[Q(z)=b_dz^d+b_{d-1}z^{d-1}+\cdots+b_2z^2+b_1z+b_0\]
in $\calg k[z]$ be the formal primitive of
\[d\!\cdot\!b_d\!\cdot\!(z-w_1')(z-w_2')\cdots(z-w_{d-1}')\]
such that $Q(0)=b_0.$ We will show that $Q$ is a polynomial with the required properties.

By the condition on $\e$ in (\ref{eq:epsilon}) and the construction of $Q$, we have that for all $1\le j\le d-1,$ the degree of $w_j'$ as a critical point of $Q$ coincides with the degree
of $w_j$ as a critical point of $P.$ Moreover,
\[\Crit(P)\cap\{\diam(y)>\e\}=\Crit(Q)\cap\{\diam(y)>\e\},\]
again by (\ref{eq:epsilon}), we have that 
\[\Crit(P)\cap\{\diam(y)>R\}=\Crit(Q)\cap\{\diam(y)>R\}.\]

It follows that $Q$ is a tame polynomial and that $\degl x=\deg_x(Q)$ for all $x\in\HK$ with $\diam(x)>R.$

By the definition of $Q$ we have that
\[b_j=\frac{b_d\!\cdot\!d}{j}\cdot\sum_{1\le i_1<i_2<\cdots<i_{d-j}\le{d-1}}w'_{i_1}w'_{i_2}\cdots w'_{i_{d-j}}\]
for $1\le j\le d-1.$

Note that if $\chara(\wtilde K)=0$ then $\abs{d}=1.$ Since $P$ is a tame polynomial it follows that if $\chara(\wtilde K)=p>0$ then $\chara(\wtilde K)$ does not divides $d,$ therefore $\abs{d}$ is also 1 in this case.

Applying Lemma~\ref{lemma:cuenta} we obtain that
\begin{align*}
\abs{b_j-a_j}&=\abs{\frac{d}{j}\cdot\sum_{1\le i_1<i_2<\cdots<i_{d-j}\le{d-1}}(a_d\cdot w_{i_1}w_{i_2}\cdots w_{i_{d-j}}-b_d\cdot w_{i_1}'w_{i_2}'\cdots w_{i_{d-j}}')}\\
&=\abs{j^{-1}}\cdot\abs{\sum_{1\le i_1<i_2<\cdots<i_{d-j}\le{d-1}}(a_d\cdot w_{i_1}w_{i_2}\cdots w_{i_{d-j}}-b_d\cdot w_{i_1}'w_{i_2}'\cdots w_{i_{d-j}}')}\\
&\le\abs{j^{-1}}\cdot\max_{1\le i_1<i_2<\cdots<i_{d-j}\le{d-1}}\abs{a_d\cdot w_{i_1}w_{i_2}\cdots w_{i_{d-j}}-b_d\cdot w_{i_1}'w_{i_2}'\cdots w_{i_{d-j}}'}\\
&\le\abs{j^{-1}}\cdot\bigg(\max_{1\le j\le{d-1}}\{\abs{a_d},\abs{w_1},\abs{w_2},\dots,\abs{w_{d-1}}\}\bigg)^{d-j}\cdot\e
\end{align*}
for all $1\le j\le d-1.$ Hence, by (\ref{eq:C}) and (\ref{eq:epsilon}), we obtain that
\begin{align}
\abs{b_j-a_j}&\le\max_{1\le j\le d-1}\Big\{\abs{j^{-1}}\Big\}\cdot\max\bigg\{1,\abs{a_d}^d,\Big(\max_{1\le j\le d-1}\{\abs{w_j}\}\Big)^d\bigg\}\cdot\e\nonumber\\
&=C\cdot\e\nonumber\\
&<\frac{R}{\max\{1,R_0^d\}}\label{eq:cer}
\end{align}
for all $0\le j\le d.$

Let $R<r<r_P$ and consider $B=B_r^+(a)$ contained in $D_0=B_{r_P}(z_0)$ such that $P(B)$ is also contained in $D_0$ and $\diam(P(B))>R.$

By (\ref{eq:r0}) and (\ref{eq:cer}), it follows that
\begin{align*}
\sup_{z\in B}\abs{(Q-P)(z)}&\le\sup_{z\in B}\max_{1\le j\le d-1}\abs{b_j-a_j}\abs z^j\\
&\le \max_{1\le j\le d-1}\abs{b_j-a_j}\cdot\max\{1,R_0^d\}\\
&<R.
\end{align*}
Therefore
\[Q(B)=P(B)+(Q-P)(B)=P(B),\]
that is, $Q(x)$ coincides with $P(x)$ for all type II points $x$ contained in 
\[\barra{D_0}\cap\{\diam(y)>R\}\cap P^{-1}(\{\diam(y)>R\}).\]

Considering decreasing sequences of type II points, we have that $Q(x)$ coincides with $P(x)$ for all $x\in \barra{D_0}\cap\{\diam(y)>R\}\cap P^{-1}(\{\diam(y)>R\}).$
\end{proof}

\begin{lemma}\label{lemma:rec_in_l}
Let $P\in K[z]$ be a tame of degree $d\ge2.$ Suppose that there exists a nonperiodic recurrent critical element $c\in\JP.$
Then $c$ belongs to $ K.$
\end{lemma}
\begin{proof}
Seeking a contradiction, suppose that there exists a recurrent and nonperiodic critical element $c\in\JP\cap\HK.$ In view of the recurrence
of $c$ we have that $P^n(c)$ belongs to the $\w$-limit of $c$ for all $n\in\N.$ Following Proposition~\ref{prop:diam} we obtain that
\[\dH(L_0,P^n(c))\le d^{d-1}\cdot\dH(L_0,c),\]
for all $n\in\N.$ In particular, $0<\inf\{\diam(P^n(c))\mid n\in\N\cup\{0\}\}.$

Let $0<R<\inf\{\diam(P^n(c))\mid n\in\N\cup\{0\}\}.$ By Proposition~\ref{prop:aprox} we have that there exists a tame polynomial $Q\in\calg k[z]$ of degree $d$ such that $Q$ coincides with $P$ in
\[\barra{D_0}\cap\{\diam(y)>R\}\cap P^{-1}(\{\diam(y)>R\}).\]
In particular, the dynamical sequence $(L_n(c))$ of $P$ is also a dynamical sequence of $Q,$ it follows that $c$ belongs to the Julia set of $Q.$ Then, $c\in\mathcal J_Q\setminus K$ is a
nonperiodic and recurrent algebraic element (since it is critical), which contradicts Proposition~\ref{prop:algnowand}. Therefore $c$ belongs to $\JPI.$
\end{proof}

\begin{proof}[Proof of Proposition~\ref{prop:critrat}]
We proceed by contradiction. Suppose that there exists a critical element $c\in\JP\cap\HK$ that is not preperiodic. According to Corollary~\ref{cor:hpcritrec} there
exists a recurrent critical element $c_1$ in $\w(c)\cap\HK$ which, by Proposition~\ref{prop:pointsperiodic}, is
not a periodic point. Applying Lemma~\ref{lemma:rec_in_l} we have that $c_1\in K,$ which is impossible.

Therefore $c$ is a preperiodic critical element. The proposition follows since periodic critical elements are
of type II by Proposition~\ref{prop:critratintro}. 
\end{proof}

\subsection{Proof of Theorem~\ref{theo:juliaberk} and Corollary~\ref{cor:nowandering}}\label{sec:proof}

In this subsection we prove a slightly stronger version of Theorem~\ref{theo:juliaberk} and we obtain some corollaries.

\begin{theorem}\label{theo:juliaberkforte}
Let $P$ be a tame polynomial with coefficients in $ K$ of degree $d\ge2.$ Then $\JP\setminus K$ is empty or
\[\JP\setminus K=\GO(x_1)\sqcup\cdots\sqcup\GO(x_m),\]
where $1\le m\le d-2$ and $x_1,\dots,x_m$ are periodic critical elements.
\end{theorem}
\begin{proof}
Consider $x\in\JP\cap\HK.$ From Corollary~\ref{cor:hpcritrec} we have that $\w(x)$ contains a nonclassical recurrent critical element. Using
Proposition~\ref{prop:critrat} and Proposition~\ref{prop:pointsperiodic} we have that $x$ is in the backward orbit of a periodic critical
element. Following Proposition~\ref{prop:finitecritical} there exist, at most, $d-2$ critical elements contained in the Julia set of $P.$ Now the
theorem follows.
\end{proof}

Corollary~\ref{cor:nowandering} follows directly from Theorem~\ref{theo:juliaberkforte} applying Proposition~\ref{prop:fatoujulia}.

If we consider a tame polynomial $P\in K[z]$ but we study the its action in the spherical completion of $ K,$ we obtain again Corollary~\ref{cor:nowandering}. 
This is not true for nontame polynomials. The example 6.3 in \cite{tesis_aste_rivera} shows that for
\[f(z)=\frac{1}{p}(z^p-z^{p^2})\in\C_p[z],\]
which is not a tame polynomial, we have that
\[J_f\sub\{x\in\ACP\mid\diam(x)=p^{-\frac{1}{p-1}}\},\]
moreover, the type II points in $J_f$ are preperiodic and the type IV points in $J_f$ are not preperiodic. Hence, if we consider the action of $f$ on the
affine line in the sense of Berkovich associated with the spherical completion of $\C_p$ we have that there exist wandering domains which are not attracted to an attracting cycle.

The following corollary is about the equilibrium measure and the entropy of $P.$ The first statement follows from the countability of $\JP\cap\HK.$ The
second statement is a direct consequence of Theorem~D in \cite{theorie_favriv}.

\begin{corollary}\label{cor:supp}
Let $P\in K[z]$ be a tame polynomial of degree $d\ge2.$ Then, the following statements hold:
\begin{enumerate}
\item The topological support of the equilibrium measure of $P$ is the classical Julia set of $P.$

\item The equilibrium measure $\rho_P$ is a measure of maximal entropy and
\[h_{\rho_P}=\htop=\log(d).\]
\end{enumerate}
\end{corollary}

Theorem~\ref{theo:juliaberkforte} and Corollary~\ref{cor:supp} are not valid for rational maps, see examples in \cite{theorie_favriv}.

It is not known if for any polynomial $P\in K[z]$ such that $\JP\cap K\ne\vacio,$ there exists a classical
repelling periodic point of $P$ in $\JP\cap K.$ In the case of tame polynomials we have the following result.
 
\begin{corollary}
Let $P\in K[z]$ be a tame polynomial of degree $d\ge2.$ Then the classical Julia set of $P$ contains a repelling periodic point.
\end{corollary}
\begin{proof}
Following Theorem~\ref{theo:juliaberkforte}, there exists $0\le m\le d-2$ such that
\[\JP\setminus K=\GO(x_1)\sqcup\cdots\sqcup\GO(x_m),\]
where $x_j$ are periodic critical elements. In particular, there exists $N\in\N,$ and a level $N$ point $L_N$ such that $L_N$ is not a critical point. If $B$ is
the ball associate to $L_N$ we have $B\sub D_0$ and $P^N(B)=D_0.$ Hence, there exists a periodic point $p\in B$ of period at most $N.$ It follows that $p$ belongs to $\JP$
and therefore $p\in K$ is a repelling periodic point.    
\end{proof}

In the case of $p$-adic polynomials, Bezivin~\cite{compacite_bezivin} proved that if there exists a classical repelling
periodic point in $\JPI$ and $\JPI$ is a compact set, then $\JP=\JPI.$ The following corollary is an analog of Proposition~A
in~\cite{compacite_bezivin}. In our case, we do not need to assume that $\JPI$ is a nonempty and compact set.

\begin{corollary}
Let $P\in K[z]$ be a tame polynomial of degree $\ge2.$ Then the following statements are equivalent:
\begin{enumerate}
\item $\JPI=\JP.$

\item All the periodic points of $P$ in $ K$ are repelling.

\item There is no critical periodic element in $\JP.$
\end{enumerate}
\end{corollary}

\subsection*{Acknowledgments}
I would like to thank Juan Rivera Letelier for suggesting to study algebraic Julia sets and for useful discussions. I am grateful to Charles Favre
and Jean-Yves Briend for their many comments and for their hospitality during my trips to France in November 2008 and in December 2009. I would like to thanks Eduardo
Friedman for useful conversation on local fields and Krasner's lemma. I am grateful for the hospitality I received from Facultad de
Ciencias de Universidad de Chile. I wish to thank especially Jan Kiwi for thoroughly revising the first draft of
this paper, several helpful suggestions and for numerous mathematical discussions.

%\bibliography{references}

\begin{thebibliography}{10}

\bibitem{equi_bakrum}
Matthew~H. Baker and Robert Rumely.
\newblock Equidistribution of small points, rational dynamics, and potential
  theory.
\newblock {\em Ann. Inst. Fourier (Grenoble)}, 56(3):625--688, 2006.

\bibitem{thesis_benedetto}
Robert~L. Benedetto.
\newblock {\em Fatou Components in $p$-adic Dynamics}.
\newblock PhD thesis, Brown University, May 1998.

\bibitem{examples_benedetto}
Robert~L. Benedetto.
\newblock Examples of wandering domains in {$p$}-adic polynomial dynamics.
\newblock {\em C. R. Math. Acad. Sci. Paris}, 335(7):615--620, 2002.

\bibitem{wanreduc_benedetto}
Robert~L. Benedetto.
\newblock Wandering domains and nontrivial reduction in non-{A}rchimedean
  dynamics.
\newblock {\em Illinois J. Math.}, 49(1):167--193 (electronic), 2005.

\bibitem{spectral_berkovich}
Vladimir~G. Berkovich.
\newblock {\em Spectral theory and analytic geometry over non-{A}rchimedean
  fields}, volume~33 of {\em Mathematical Surveys and Monographs}.
\newblock American Mathematical Society, Providence, RI, 1990.

\bibitem{compacite_bezivin}
Jean-Paul B{\'e}zivin.
\newblock Sur la compacit{\'e} des ensembles de {J}ulia des polyn{\^o}mes
  {$p$}-adiques.
\newblock {\em Math. Z.}, 246(1-2):273--289, 2004.

\bibitem{local_cassels}
J.~W.~S. Cassels.
\newblock {\em Local fields}, volume~3 of {\em London Mathematical Society
  Student Texts}.
\newblock Cambridge University Press, Cambridge, 1986.

\bibitem{measure_chamber}
Antoine Chambert-Loir.
\newblock Mesures et {\'e}quidistribution sur les espaces de {B}erkovich.
\newblock {\em J. Reine Angew. Math.}, 595:215--235, 2006.

\bibitem{valuative_favjon}
Charles Favre and Mattias Jonsson.
\newblock {\em The valuative tree}, volume 1853 of {\em Lecture Notes in
  Mathematics}.
\newblock Springer-Verlag, Berlin, 2004.

\bibitem{equi_favriv}
Charles Favre and Juan Rivera-Letelier.
\newblock {\'E}quidistribution quantitative des points de petite hauteur sur la
  droite projective.
\newblock {\em Math. Ann.}, 335(2):311--361, 2006.

\bibitem{theorie_favriv}
Charles Favre and Juan Rivera-Letelier.
\newblock Th{\'e}orie ergodique des fractions rationnelles sur un corps
  ultram{\'e}trique.
\newblock {\em Proc. Lond. Math. Soc. (3)}, 100(1):116--154, 2010.

\bibitem{padic_gouvea}
Fernando~Q. Gouv{\^e}a.
\newblock {\em {$p$}-adic numbers}.
\newblock Universitext. Springer-Verlag, Berlin, 1993.
\newblock An introduction.

\bibitem{weak_hsia}
Liang-Chung Hsia.
\newblock A weak {N}{\'e}ron model with applications to {$p$}-adic dynamical
  systems.
\newblock {\em Compositio Math.}, 100(3):277--304, 1996.

\bibitem{closure_hsia}
Liang-Chung Hsia.
\newblock Closure of periodic points over a non-{A}rchimedean field.
\newblock {\em J. London Math. Soc. (2)}, 62(3):685--700, 2000.

\bibitem{pui_kiwi}
Jan Kiwi.
\newblock Puiseux series polynomial dynamics and iteration of complex cubic
  polynomials.
\newblock {\em Ann. Inst. Fourier (Grenoble)}, 56(5):1337--1404, 2006.

\bibitem{branner_qiuyin}
WeiYuan Qiu and YongCheng Yin.
\newblock Proof of the {B}ranner-{H}ubbard conjecture on {C}antor {J}ulia sets.
\newblock {\em Sci. China Ser. A}, 52(1):45--65, 2009.

\bibitem{tesis_rivera}
Juan Rivera-Letelier.
\newblock {\em Dynamique de fractions rationnelles sur des corps locaux}.
\newblock PhD thesis, U. de Paris-Sud, Orsay, 2000.

\bibitem{tesis_aste_rivera}
Juan Rivera-Letelier.
\newblock Dynamique des fonctions rationnelles sur des corps locaux.
\newblock {\em Ast{\'e}risque}, (287):xv, 147--230, 2003.
\newblock Geometric methods in dynamics. II.

\bibitem{espace_rivera}
Juan Rivera-Letelier.
\newblock Espace hyperbolique {$p$}-adique et dynamique des fonctions
  rationnelles.
\newblock {\em Compositio Math.}, 138(2):199--231, 2003.

\bibitem{points_rivera}
Juan Rivera-Letelier.
\newblock Points p{\'e}riodiques des fonctions rationnelles dans l'espace
  hyperbolique {$p$}-adique.
\newblock {\em Comment. Math. Helv.}, 80(3):593--629, 2005.

\bibitem{notes_rivera}
Juan Rivera-Letelier.
\newblock Notes sur la droite projective de {B}erkovich.
\newblock \url{www.arxiv.org/pdf/math/0605676}, 2006.

\bibitem{quasiconformalI_sullivan}
Dennis Sullivan.
\newblock Quasiconformal homeomorphisms and dynamics. {I}. {S}olution of the
  {F}atou-{J}ulia problem on wandering domains.
\newblock {\em Ann. of Math. (2)}, 122(3):401--418, 1985.

\bibitem{algebraic_walker}
Robert~J. Walker.
\newblock {\em Algebraic curves}.
\newblock Springer-Verlag, New York, 1978.
\newblock Reprint of the 1950 edition.

\end{thebibliography}
%\bibliographystyle{plain}

\end{document}